\documentclass[11pt]{amsart}

\newtheorem{proposition}{Proposition}[section]
\newtheorem{lemma}[proposition]{Lemma}

\newtheorem{theorem}[proposition]{Theorem}

\theoremstyle{definition}

\theoremstyle{remark}
\newtheorem{remark}[proposition]{Remark}

\newtheorem{example}[proposition]{Example}

\newcommand{\thlabel}[1]{\label{th:#1}}
\newcommand{\thref}[1]{Theorem~\ref{th:#1}}
\newcommand{\selabel}[1]{\label{se:#1}}
\newcommand{\seref}[1]{Section~\ref{se:#1}}
\newcommand{\lelabel}[1]{\label{le:#1}}
\newcommand{\leref}[1]{Lemma~\ref{le:#1}}
\newcommand{\prlabel}[1]{\label{pr:#1}}
\newcommand{\prref}[1]{Proposition~\ref{pr:#1}}
\newcommand{\colabel}[1]{\label{co:#1}}

\newcommand{\relabel}[1]{\label{re:#1}}

\newcommand{\exlabel}[1]{\label{ex:#1}}
\newcommand{\exref}[1]{Example~\ref{ex:#1}}

\newcommand{\eqlabel}[1]{\label{eq:#1}}
\newcommand{\equref}[1]{(\ref{eq:#1})}

\def\equal#1{\smash{\mathop{=}\limits^{#1}}}

\newcommand{\can}{{\rm can}}
\newcommand{\Hom}{{\rm Hom}}
\newcommand{\End}{{\rm End}}

\newcommand{\im}{{\rm Im}\,}

\def\lan{\langle}
\def\ran{\rangle}
\def\ot{\otimes}

\def\leftact{\hbox{$\rightharpoonup$}}

\newcommand{\Cc}{\mathcal{C}}

\newcommand{\Mm}{\mathcal{M}}

\newcommand{\Rr}{\mathcal{R}}

\def\*C{{}^*\hspace*{-1pt}{\Cc}}

\def\text#1{{\rm {\rm #1}}}

\def\ol{\overline}
\def\ul{\underline}

\begin{document}
\title[Partial (co)actions of Hopf algebras]{Partial (co)actions of Hopf algebras and partial Hopf-Galois theory }
\author{S. Caenepeel}
\address{Faculty of Engineering, 
Vrije Universiteit Brussel, B-1050 Brussels, Belgium}
\email{scaenepe@vub.ac.be}
\urladdr{http://homepages.vub.ac.be/\~{}scaenepe/}
\author{K. Janssen}
\address{Faculty of Engineering,
Vrije Universiteit Brussel, B-1050 Brussels, Belgium}
\email{krjansse@vub.ac.be}
\urladdr{http://homepages.vub.ac.be/\~{}krjansse/}
\thanks{This research was supported by the research project G.0622.06 ``Deformation quantization methods
for algebras and categories with applications to quantum mechanics" from
FWO-Vlaanderen.}



\begin{abstract}
We introduce partial (co)actions of a Hopf algebra $H$ on an algebra. To this end,
we introduce first the notion of lax coring, generalizing Wisbauer's notion of
weak coring. We also have the dual notion of lax ring. 
Several duality results are given, and we develop
Galois theory for partial $H$-comodule algebras.
\end{abstract}
\maketitle
\section*{Introduction}
Partial group actions were considered first by Exel \cite{Exel}, in the context of
operator algebras. A treatment from a purely algebraic point of view was given recently
in \cite{DEP,DE,DFP,DZ}. In particular, Galois theory over commutative rings can
be generalized to partial group actions, see \cite{DFP} (at least under the additional
assumption that the associated ideals are generated by idempotents).\\
The following questions arise naturally: can we develop a theory of partial (co)actions
of Hopf algebras? Is it possible to generalize Hopf-Galois theory to the partial situation?
The aim of this paper is to give a positive answer to these questions, with one important
restriction: our approach only leads to a generalization of partial group actions, with
associated ideals generated by central idempotents. \\
Partial group actions were studied from the point of view of corings by the first author
and De Groot in \cite{CaenepeelDG05}. Namely, a partial group action in
the sense of \cite{DFP} gives rise to a coring. The Galois theory of \cite{DFP}
can then be considered as a special case of the Galois theory of corings
(see \cite{Br3,BrW,Caenepeel03,W3}). There is a remarkable analogy with the
Galois theory that can be developed for weak Hopf algebras (see \cite{CDG2}):
in both cases, the associated coring is a direct factor of the tensor product of
the Galois extension $A$, and a coalgebra. In the partial group action case, the
coalgebra is the dual of the group algebra, in the other case it is the weak Hopf algebra
that we started with. The right $A$-module structure of the coring is induced by a
kind of entwining map. In the weak Hopf algebra case, it is a weak entwining map,
as introduced in \cite{CaenepeelDG00}. The map in the partial group action case,
however, does not satisfy the axioms of a weak entwining structure.\\
Wisbauer \cite{Wisbauer} introduces weak entwining structures from the point of
view of {\sl weak} corings; these are corings with a bimodule structure that is not
necessarily unital. If $\Cc$ is a left-unital weak $A$-coring, then $\Cc 1_A$ is
an $A$-coring that is a direct summand of $\Cc$. Weak entwining structures are
then in one-to-one correspondence with left-unital weak $A$-coring structures
on $A\ot C$, where $A$ is an algebra, and $C$ is a coalgebra.\\
If a finite group $G$ acts partially on an algebra $A$, then we can define a
left-unital $A$-bimodule structure on $A\ot (kG)^*$, such that $(A\ot (kG)^*)1_A$
is an $A$-coring, and a direct factor of $A\ot (kG)^*$. But $A\ot (kG)^*$ does not
satisfy Wisbauer's axioms of a weak coring. This observation has lead us to the
introduction of {\sl lax corings}. The counit property of a lax coring is weaker than
that of a weak coring, but it is still designed in such a way that $\Cc 1_A$ is a
coring.\\
Now let $H$ be a Hopf algebra, and consider a map $\rho: A\to A\ot H$.
Our next step is then to examine lax coring structures on $A\ot H$. A subtlety that
appears is that we have two possible choices for the counit: we can consider
$A\ot \epsilon$ and $(A\ot \epsilon)\circ \pi$, where $\epsilon$ is the counit on $H$,
and $\pi$ is the projection of $A\ot H$ onto $(A\ot H)1_A$. This leads to the
introduction of a right partial (resp. lax) $H$-comodule algebra $A$. The notion of lax comodule algebra is
the most general, and includes partial and weak comodule algebras as special
cases. If $A$ is at the same time a partial and weak comodule algebra,
then it is a comodule algebra.\\
We have a dual theory: we can introduce lax $A$-rings, and we then obtain the definition of partial (resp. lax)
$H$-module algebra. In the case where $H$ is a group algebra, we recover the
definition of partial group action. We also discuss duality results. For example, if
$H$ is a finitely generated projective bialgebra, then we have a bijective correspondence
between right lax (resp. partial) $H$-comodule algebra structures on $A$ and left lax (resp. partial) $H^{*{\rm cop}}$-module algebra structures on $A^{\rm op}$ (see \thref{4.8}).
In the final \seref{7}, we applied the theory of Galois corings to corings arising
from partial comodule algebras.

\section{Lax rings and corings}\selabel{1}
Let $A$ be a ring with unit. $A$-modules will not necessarily be unital.

\begin{proposition}\prlabel{1.1}
Let $P$ be a unital left $A$-module. There is a bijective correspondence between
(non-unital) right $A$-module structures on $P$ making $P$ an $A$-bimodule and
unital right $A$-module structures on left $A$-linear direct factors $\ul{P}$
of $P$, making $\ul{P}$ a unital $A$-bimodule.
\end{proposition}

\begin{proof}
For an $A$-bimodule $P$, the map $\pi: P\to P$, $\pi(p)=p 1_A$
is a left $A$-linear projection. The right $A$-action on $P$ restricts to a 
unital right $A$-action on $\ul{P}=\im (\pi)$.
Conversely, let $\pi: P\to \ul{P}$ be a left $A$-linear projection, and let $\ul{P}$
be a unital $A$-bimodule. We extend the right $A$-action from $\ul{P}$ to $P$
as follows:
$pa=\pi(p)a\in \ul{P}$.
This action is associative, since
$(pa)b=(\pi(p)a)b=\pi(p)(ab)=p(ab)$.
\end{proof}

We observe that $\pi$ is then also right $A$-linear, so $\ul{P}$ is an $A$-bimodule
direct factor of $P$. The inclusion $\iota: \ul{P}\to P$ is a right inverse of $\pi$.

Recall that an $A$-coring $(\Cc,\Delta,\varepsilon)$ is a coalgebra in the monoidal category ${}_A\Mm_A$ of unital $A$-bimodules. This means that $\Cc$ is a unital $A$-bimodule, and
that $\Delta: \Cc\to \Cc\ot_A\Cc$ and $\varepsilon: \Cc\to A$ are $A$-bimodule maps
such that
\begin{eqnarray}
&&(\Delta\ot_A\Cc)\circ\Delta=(\Cc\ot_A\Delta)\circ\Delta;\eqlabel{1.2.1}\\
&&(\varepsilon\ot_A\Cc)\circ\Delta=(\Cc\ot_A\varepsilon)\circ\Delta= \Cc.\eqlabel{1.2.2}
\end{eqnarray}
We will use the Sweedler-Heyneman notation
$\Delta(c)=c_{(1)}\ot_A c_{(2)}$,
where summation is implicitly understood.

Now take a left unital $A$-bimodule $\Cc$ and two $A$-bimodule maps
$\Delta: \Cc\to \Cc\ot_A\Cc$ and $\varepsilon: \Cc\to A$ satisfying \equref{1.2.1}.
We consider the projection $\pi: \Cc\to\ul{\Cc}=\Cc 1_A$ and its right inverse $\iota$.
$\Delta$ restricts to a map $\ul{\Delta}: \ul{\Cc}\to \ul{\Cc}\ot_A\ul{\Cc}$ since
$\Delta(c1_A)=\Delta(c)1_A=c_{(1)}\ot_A 1_Ac_{(2)}1_A=c_{(1)}1_A\ot_A c_{(2)}1_A
\in \ul{\Cc}\ot_A\ul{\Cc}$, for all $c\in \Cc$.
$\varepsilon\circ\iota$
is then the restriction of $\varepsilon$ to $\ul{\Cc}$.\\
We call $(\Cc,\Delta,\varepsilon)$ a left unital {\sl lax} (resp. {\sl weak}) $A$-coring if
\equref{1.2.3} (resp. \equref{1.2.4}) holds for all $c\in \ul{\Cc}$ (resp. $c\in \Cc$).
\begin{eqnarray}
c&=&{\varepsilon}(c_{(1)})c_{(2)}=c_{(1)}{\varepsilon}(c_{(2)});\eqlabel{1.2.3}\\
c1_A&=&\varepsilon(c_{(1)})c_{(2)}=c_{(1)}\varepsilon(c_{(2)}).\eqlabel{1.2.4}
\end{eqnarray}
Weak corings were introduced in \cite{Wisbauer}. $(\Cc,\Delta,\varepsilon)$ is
a left unital lax $A$-coring if and only if 
$(\ul{\Cc},\ul{\Delta},{\varepsilon}\circ \iota)$ is an $A$-coring. Clearly weak corings
are lax.

Recall that an $A$-ring $(\Rr,\mu,\eta)$ is an algebra in the category of unital
$A$-bimodules. This means that $\mu: \Rr\ot_A\Rr\to \Rr$ and $\eta: A\to \Rr$
are $A$-bimodule maps such that
\begin{eqnarray}\eqlabel{1.3.1}
&&\mu\circ (\mu\ot_A \Rr)=\mu\circ (\Rr\ot_A \mu)\\
&&\eqlabel{1.3.2}
\mu\circ (\eta\ot_A \Rr)=\mu\circ (\Rr\ot_A \eta)=\Rr.
\end{eqnarray}
Then $\Rr$ is a ring with unit $\eta(1_A)$, and $\eta: A\to \Rr$ is a ring morphism.
It follows from \equref{1.3.2} that the $A$-bimodule structure on $\Rr$ is induced by $\eta$.
So an $A$-ring is a ring $\Rr$ together with a
ring morphism $\eta: A\to \Rr$.\\
Let $\Rr$ be a left unital $A$-bimodule, and consider the projection 
$\pi: \Rr\to\ul{\Rr}=\Rr 1_A$, and
an $A$-bimodule map $\mu: \Rr\ot_A\Rr\to \Rr$ satisfying \equref{1.3.1}.
$\mu$ restricts to a map
$\ul{\mu}: \ul{\Rr}\ot_A\ul{\Rr}\to \ul{\Rr}$, since
$\mu(r1_A\ot_A s1_A)=\mu(r1_A\ot_A s)1_A\in \ul{\Rr}$, for all $r,s\in \Rr$.
We will write $\mu(r\ot_A s)=rs$, as usual.
Let $\eta: A\to \Rr$ be an $A$-bimodule map, and write
$\eta(1_A)=1_{\Rr}$. Then $\pi(1_\Rr)=1_\Rr1_A=\eta(1_A)1_A=\eta(1_A)=1_\Rr$,
so that $\pi \circ \eta$ is the corestriction of $\eta$ to $\ul{\Rr}$.\\
$(\Rr,\mu,\eta)$ is called a left unital {\sl lax} (resp. {\sl weak}) $A$-ring if
\equref{1.4.1} (resp. \equref{1.4.2}) holds for all $r\in \ul{\Rr}$ (resp. for all $r\in {\Rr}$).
\begin{eqnarray}
r&=&{1}_\Rr r=r{1}_\Rr;\eqlabel{1.4.1}\\
r1_A&=&{1}_\Rr r=r1_\Rr.\eqlabel{1.4.2}
\end{eqnarray}
$(\Rr,\mu,\eta)$ is a left unital lax $A$-ring if and only if $(\ul{\Rr},\ul{\mu}, {\pi \circ \eta})$
is an $A$-ring.\\
Right unital lax and weak $A$-rings are introduced in a similar way. Let 
 $\Rr$ be a right unital $A$-bimodule and $\ul{\Rr}=1_A \Rr$. Consider
an $A$-bimodule map $\mu: \Rr\ot_A\Rr\to \Rr$ satisfying \equref{1.3.1}.
$\mu$ restricts to $\ul{\mu}: \ul{\Rr}\ot_A\ul{\Rr}\to \ul{\Rr}$. 
$\eta: A\to \Rr$ corestricts to the map $\pi \circ \eta: A\to \ul{\Rr}$.
$(\Rr,\mu,\eta)$ is a right unital {\sl lax} (resp. {\sl weak}) $A$-ring if
\equref{1.4.1} is fulfilled for all $r\in \ul{\Rr}$ (resp. $1_A r={1}_\Rr r=r1_\Rr$
for all $r\in {\Rr}$).

\subsubsection*{Duality}
Let $\Cc$ be a left unital $A$-bimodule. Then ${}^*\Cc={}_A\Hom(\Cc,A)$ is a right
unital $A$-bimodule, with $A$-action
$(afb)(c)=f(ca)b,$
for all $a,b\in A$, $c\in \Cc$ and $f\in {}^*\Cc$. If $\Delta: \Cc\to \Cc\ot_A\Cc$
is a coassociative $A$-bimodule map, then $\mu: {}^*\Cc\ot_A{}^*\Cc\to {}^*\Cc$,
$\mu(f\ot_A g)=f \# g$ given by
$(f\# g)(c)=g(c_{(1)}f(c_{(2)}))$,
for all $c\in \Cc$, is an associative $A$-bimodule map.
Let $\varepsilon: \Cc\to A$ be an $A$-bimodule map. For all $a\in A$ and $c\in \Cc$,
we have that
$(a\varepsilon)(c)=(\varepsilon a)(c)=\varepsilon(c)a$,
so $\eta: A\to {}^*\Cc$, $\eta(a)=a\varepsilon=\varepsilon a$, is an $A$-bimodule map.
For all $f\in {}^*\Cc$ and $c\in \Cc$, we compute that
$(\varepsilon\# f)(c)=f(c_{(1)}\varepsilon(c_{(2)}))$
and
$(f\#\varepsilon)(c)=
f(\varepsilon(c_{(1)})c_{(2)}).$

\begin{proposition}\prlabel{1.5}
If $(\Cc,\Delta,\varepsilon)$ is a left unital weak (resp. lax) $A$-coring, then
$({}^*\Cc,\mu,\eta)$ is a right unital weak (resp. lax) $A$-ring. The $A$-rings $1_A{}^*\Cc$ and ${}^*(\Cc1_A)$ are isomorphic.
\end{proposition}

\begin{proof}
Let $(\Cc,\Delta,\varepsilon)$ be a left unital weak $A$-coring. 
Then $({}^*\Cc,\mu,\eta)$ is a right unital weak $A$-ring since
$(\varepsilon\# f)(c)= f(c_{(1)}\varepsilon(c_{(2)}))=
f(c1_A)=(1_A f)(c)$ and
$(f\#\varepsilon)(c) = f(\varepsilon(c_{(1)})c_{(2)})
= f(c1_A)=(1_A f)(c)$,
for all $f\in {}^*\Cc$ and $c\in \Cc$.\\
Now assume that $(\Cc,\Delta,\varepsilon)$ is a left unital lax $A$-coring. 
For all $\ul{f}=1_A f\in 1_A{}^*\Cc=\ul{{}^*\Cc}$ and $c\in \Cc$, we have
\begin{eqnarray*}
&&\hspace*{-2cm}
(\varepsilon \# \ul{f})(c)= (\varepsilon \# 1_A f)(c) = (1_A f)(c_{(1)}\varepsilon(c_{(2)}))\\
&=& f(c_{(1)}\varepsilon(c_{(2)}) 1_A) = f(c_{(1)}\varepsilon(c_{(2)}))\equal{\equref{1.2.3}} f(c 1_A)= (1_A f)(c)\\
&&\hspace{-2cm}
(\ul{f}\# \varepsilon)(c) = \varepsilon(c_{(1)}(1_A f)(c_{(2)}))= \varepsilon(c_{(1)}f(c_{(2)}1_A))\\
&=& \varepsilon(c_{(1)})f(c_{(2)}1_A) = f(\varepsilon(c_{(1)})c_{(2)}1_A) \equal{\equref{1.2.3}} f(c 1_A) = (1_A f)(c).
\end{eqnarray*}
So $\varepsilon \# \ul{f}=\ul{f}\#\varepsilon =\ul{f}$, and \equref{1.4.1} holds, and it follows that
$({}^*\Cc,\mu,\eta)$ is a right unital lax $A$-ring.\\
To prove the final statement, we observe first that $f\in 1_A{}^*\Cc$ if and only if $1_A\cdot f=f$,
or $f(c 1_A)=f(c)$, for all $c\in \Cc$. Consider the maps
$\alpha: 1_A{}^*\Cc\to {}^*(\Cc1_A),~~\alpha(f)=f_{|\Cc1_A},$
and $\beta: {}^*(\Cc1_A) \to 1_A{}^*\Cc,~~\beta(g)(c)=g(c1_A)$. It is easily verified that $\beta(g)\in 1_A{}^*\Cc$, and
that $\alpha$ and $\beta$ are inverses.
\end{proof}

\section{Partial comodule algebras}\selabel{2}
Let $k$ be a commutative ring and $A,H$ be two $k$-algebras. 
$A\ot H$ is a unital left $A\ot H$-module via the multiplication on $A\ot H$, and a left unital $A$-module
via restriction of scalars.
For a
$k$-linear (not necessarily coassociative) map $\rho: A\to A\ot H$, we adopt the notation
$$\rho(a)=a_{[0]}\ot a_{[1]}=a_{[0']}\ot a_{[1']},~ (\rho\ot H)(\rho(a))=\rho^2(a)=a_{[0]}\ot a_{[1]}\ot a_{[2]},$$
for all $a\in A$. Summation is implicitly understood.

\begin{lemma}\lelabel{2.1}
There is a bijective correspondence between
 (non-unital) right $A$-actions on $A\ot H$, compatible with the left $A\ot H$-action,
 and $k$-linear maps $\rho: A\to A\ot H$ satisfying 
\begin{equation}\eqlabel{2.1.1}
\rho(ab)=\rho(a)\rho(b)=a_{[0]}b_{[0]}\ot a_{[1]}b_{[1]},
\end{equation}
for all $a,b\in A$.
Then $\ul{A\ot H}=(A\ot H)1_A$ is a unital $A$-bimodule, and we have a projection
$$\pi: A\ot H\to A\ot H,~~\pi(a\ot h)=(a\ot h)1_A=a1_{[0]}\ot h1_{[1]}.$$
$A\ot H$ is right $A$-unital if
\begin{equation}\eqlabel{2.1.2}
\rho(1_A)=1_A\ot 1_H.
\end{equation}
\end{lemma}

\begin{proof}
Given a right $A$-action, we define $\rho: A\to A\ot H$ by the formula
$\rho(a)=(1_A\ot 1_H)a$.
\equref{2.1.1} follows from the associativity of the right $A$-action. Conversely,
given $\rho$, we define a right $A$-action by $(a\ot h)b=ab_{[0]}\ot hb_{[1]}$.
\end{proof}

Now let $(H,\delta,\epsilon)$ be a $k$-bialgebra, and $\rho: A\to A\ot H$
a $k$-linear map satisfying \equref{2.1.1}. Consider the left $A$-linear maps
\begin{eqnarray*}
&&\hspace*{-1cm}\Delta=(\pi\ot H)\circ(A\ot \delta): A\ot H\to A\ot H\ot_AA\ot H\cong \ul{A\ot H}\ot H,\\
&&
\Delta(a\ot h)=(a\ot h_{(1)})\ot_A (1_A\ot h_{(2)})\cong a1_{[0]} \ot h_{(1)}1_{[1]}\ot h_{(2)},\\
&&\hspace*{-1cm}{\varepsilon}=A\ot \epsilon: {A\ot H}\to A,~~
{\varepsilon}(a\ot h)=\epsilon(h)a,\\
&&\hspace*{-1cm}\ul{\varepsilon}=(A\ot \epsilon)\circ\pi: {A\ot H}\to A,~~
\ul{\varepsilon}(a\ot h)=\epsilon(h1_{[1]})a1_{[0]}=\epsilon(h)\epsilon(1_{[1]})a1_{[0]}.
\end{eqnarray*}

We will now investigate when $({A\ot H},\Delta,\varepsilon)$ and
$({A\ot H},\Delta,\ul{\varepsilon})$ are left unital weak, resp. lax $A$-corings.
Then $\Delta$ and $\varepsilon$ (or $\ul{\varepsilon}$) have to be right $A$-linear.

\begin{lemma}\lelabel{2.2}
Let $A$ be a $k$-algebra, $H$ a $k$-bialgebra, and $\rho: A\to A\ot H$
a $k$-linear map satisfying \equref{2.1.1}.\\
1) $\Delta$ is right $A$-linear if and only if, for all $a\in A$,
\begin{equation}\eqlabel{2.2.1}
\rho^2(a)=a_{[0]}1_{[0]}\ot a_{[1](1)}1_{[1]}\ot a_{[1](2)}.
\end{equation}
2) $\ul{\varepsilon}$ is right $A$-linear if and only if, for all $a\in A$,
\begin{equation}\eqlabel{2.2.2}
\epsilon(a_{[1]})a_{[0]}=\epsilon(1_{[1]})1_{[0]}a.
\end{equation}
3) $\varepsilon$ is right $A$-linear if and only if $\varepsilon=\ul{\varepsilon}$,
that is, for all $a\in A$,
\begin{equation}\eqlabel{2.2.3}
\epsilon(a_{[1]})a_{[0]}=a.
\end{equation}
\end{lemma}

\begin{proof}
1) follows immediately from the following observation:
\begin{eqnarray*}
&&\hspace*{-1cm}
\Delta(1_A\ot h)a=(1_A\ot h_{(1)})\ot_A (a_{[0]} \ot h_{(2)}a_{[1]})\cong
a_{[0]}\ot h_{(1)}a_{[1]}\ot h_{(2)}a_{[2]};\\
&&\hspace*{-1cm}
\Delta((1_A\ot h)a)=\Delta(a_{[0]}\ot ha_{[1]})=(a_{[0]}\ot (ha_{[1]})_{(1)})\ot_A (1_A\ot (ha_{[1]})_{(2)})\\
&=&(a_{[0]}\ot h_{(1)}a_{[1](1)})\ot_A (1_A\ot h_{(2)}a_{[1](2)})\\
&\cong& a_{[0]}1_{[0]}\ot h_{(1)}a_{[1](1)}1_{[1]}\ot h_{(2)}a_{[1](2)}.
\end{eqnarray*}
2) $\ul{\varepsilon}$ is right $A$-linear if and only if $\ul{\varepsilon}(1_A\ot h)a=\epsilon(h)\epsilon(1_{[1]})1_{[0]}a$ equals $\ul{\varepsilon}((1_A\ot h)a)=\ul{\varepsilon}(a_{[0]}\ot ha_{[1]})=\epsilon(ha_{[1]}1_{[1]})a_{[0]}1_{[0]}=\epsilon(ha_{[1]})a_{[0]}=\epsilon(h)\epsilon(a_{[1]})a_{[0]}$,
for all $a\in A$ and $h\in H$, and this is equivalent to \equref{2.2.2}.\\
3) $\varepsilon$ is right $A$-linear if and only if $\varepsilon(1_A\ot h)a=\epsilon(h)a$ equals $\varepsilon((1_A\ot h)a)=\epsilon(ha_{[1]})a_{[0]}=\epsilon(h)\epsilon(a_{[1]})a_{[0]}$, and this is equivalent to \equref{2.2.3}.
\end{proof}

If $\Delta$ is right $A$-linear, then $\Delta$ restricts to a map
$\ul{\Delta}: \ul{A\ot H}\to \ul{A\ot H}\ot_A \ul{A\ot H}$.

\begin{proposition}\prlabel{2.3}
Let $A$ be a $k$-algebra, $H$ a $k$-bialgebra, and $\rho:A\to A\ot H$
a $k$-linear map. We call $A$ a (right) weak $H$-comodule algebra if
the following equivalent conditions hold:\\
1) $(A\ot H,\Delta,\ul{\varepsilon})$ is a left unital weak $A$-coring;\\
2) the conditions (\ref{eq:2.1.1},\ref{eq:2.2.1},\ref{eq:2.2.2},\ref{eq:2.3.1}) 
are satisfied, for all $a,b\in A$;\\
3) the conditions (\ref{eq:2.1.1},\ref{eq:2.2.2},\ref{eq:2.3.1},\ref{eq:2.3.2})
are satisfied for all $a,b\in A$.
\begin{eqnarray}
\rho(1_A)&=&\epsilon(1_{[1]})1_{[0]}\ot 1_H;\eqlabel{2.3.1}\\
\rho^2(a)&=&a_{[0]}\ot \delta(a_{[1]}).\eqlabel{2.3.2}
\end{eqnarray}
\end{proposition}

\begin{proof}
If $(A\ot H,\Delta,\ul{\varepsilon})$ is a left unital weak $A$-coring, it follows that
(\ref{eq:2.1.1},\ref{eq:2.2.1},\ref{eq:2.2.2}) hold, by Lemmas \ref{le:2.1}
and \ref{le:2.2}. The left counit property $(\ul{\varepsilon}\ot_A (A\ot H))\circ\Delta= \pi$ (cf. \equref{1.2.4})
then holds if and only if
$
\ul{\varepsilon}(a\ot h_{(1)})(1_A\ot h_{(2)})=\epsilon(h_{(1)}1_{[1]})a1_{[0]} \ot h_{(2)}=\epsilon(1_{[1]})a1_{[0]}\ot h$
equals
$(a\ot h)1_A=a1_{[0]}\ot h1_{[1]},$
for all $a\in A$ and $h\in H$, if and only if \equref{2.3.1} holds.
This proves that $\ul{1)\Rightarrow 2)}$. It also proves that $\ul{2)\Rightarrow 1)}$,
if we can show that the
right counit condition is satisfied. Indeed,
\begin{eqnarray*}
&&\hspace*{-15mm}
(a\ot h_{(1)})\ul{\varepsilon}(1_A\ot h_{(2)})=\epsilon(h_{(2)}1_{[1]})(a\ot h_{(1)})1_{[0]}\\
&=& \epsilon(h_{(2)}1_{[2]})a1_{[0]}\ot h_{(1)}1_{[1]}
\equal{\equref{2.2.1}}\epsilon(h_{(2)}1_{[1](2)})a1_{[0]}1_{[0']}\ot h_{(1)}1_{[1](1)}1_{[1']}\\
&=& a1_{[0]}1_{[0']}\ot h1_{[1]}1_{[1']}
\equal{\equref{2.1.1}}a1_{[0]}\ot h1_{[1]}=(a\ot h)1_A.
\end{eqnarray*}
$\ul{2)\Leftrightarrow 3)}$. We will prove that the right hand sides of the formulas 
\equref{2.2.1} and \equref{2.3.2} are equal if $\rho$ satisfies \equref{2.1.1} and \equref{2.3.1}.
Indeed,
\begin{eqnarray*}
&&\hspace*{-2cm}
a_{[0]}1_{[0]}\ot a_{[1](1)}1_{[1]}\ot a_{[1](2)}\equal{ \equref{2.3.1}}
a_{[0]}\epsilon(1_{[1]})1_{[0]}\ot a_{[1](1)}\ot a_{[1](2)}\\
&\equal{\equref{2.3.1}}&
a_{[0]}1_{[0]}\ot (a_{[1]}1_{[1]})_{(1)}\ot (a_{[1]}1_{[1]})_{(2)}\\
&\equal{ \equref{2.1.1}}& a_{[0]}\ot a_{[1](1)} \ot a_{[1](2)}=a_{[0]}\ot \delta(a_{[1]}).
\end{eqnarray*}
\end{proof}

\begin{proposition}\prlabel{2.4}
Let $A$ be a $k$-algebra, $H$ a $k$-bialgebra, and $\rho: A\to A\ot H$
a $k$-linear map. The following assertions are equivalent:\\
1) $(A\ot H,\Delta,{\varepsilon})$ is an $A$-coring;\\
2) $(A\ot H,\Delta,{\varepsilon})$ is a left unital weak $A$-coring;\\
3) $A$ is a right $H$-comodule algebra; this means that the conditions
 (\ref{eq:2.1.1},\ref{eq:2.1.2},\ref{eq:2.2.3}) and \equref{2.3.2} are satisfied, for all $a,b\in A$.
\end{proposition}

\begin{proof}
$\ul {2)\Rightarrow 3)}$. 
It follows from Lemmas \ref{le:2.1} and \ref{le:2.2} that 
(\ref{eq:2.1.1},\ref{eq:2.2.1},\ref{eq:2.2.3}) hold. Using \equref{1.2.4}, we find
$\rho(1_A)=(1_A\ot1_H)1_A=\varepsilon(1_A\ot 1_H)(1_A\ot 1_H)=\epsilon(1_H)1_A\ot 1_H=1_A\ot 1_H,$
so \equref{2.1.2} holds. \equref{2.3.2} follows easily from
\equref{2.1.2} and \equref{2.2.1}.
$\ul {3)\Rightarrow 1)}$ is well-known (see e.g. \cite{Br3}). $\ul {1)\Rightarrow 2)}$ is trivial.
\end{proof}

\begin{proposition}\prlabel{2.5}
Let $A$ be a $k$-algebra, $H$ a $k$-bialgebra, and $\rho: A\to A\ot H$
a $k$-linear map. We call $A$ a (right) lax $H$-comodule algebra if the following 
equivalent assertions are satisfied:\\
1) $(A\ot H,\Delta,\ul{\varepsilon})$ is a left unital lax $A$-coring;\\
2) the conditions (\ref{eq:2.1.1},\ref{eq:2.2.1},\ref{eq:2.2.2},\ref{eq:2.5.1}) are satisfied, for all $a,b\in A$;\\
3) the conditions (\ref{eq:2.1.1},\ref{eq:2.2.1},\ref{eq:2.2.2},\ref{eq:2.5.2}) are satisfied, for all $a,b\in A$;\\
4) the conditions (\ref{eq:2.1.1},\ref{eq:2.2.1},\ref{eq:2.2.2},\ref{eq:2.5.3}) are satisfied, for all $a,b\in A$.
\begin{eqnarray}
\rho(1_A)&=&\epsilon(1_{[1]})1_{[0]}\ot 1_{[2]}\eqlabel{2.5.1};\\
\rho(1_A)&=&\epsilon(1_{[1']})1_{[0']}1_{[0]}\ot 1_{[1]}\eqlabel{2.5.2};\\
\rho(1_A)&=&\epsilon(1_{[1']})1_{[0]}1_{[0']}\ot 1_{[1]}\eqlabel{2.5.3}.
\end{eqnarray}
\end{proposition}

\begin{proof}
$\ul{1)\Rightarrow 2)}$. It follows from Lemmas \ref{le:2.1} and \ref{le:2.2}
that (\ref{eq:2.1.1},\ref{eq:2.2.1},\ref{eq:2.2.2})  hold.
We have that
\begin{eqnarray*}
&&\hspace*{-2cm}
\ul{\Delta}((1_A\ot 1_H)1_A)=\ul{\Delta}(1_{[0]}\ot 1_{[1]})=(1_{[0]}\ot 1_{[1](1)})\ot_A (1_A\ot 1_{[1](2)})\\
&=& (1_{[0]} 1_{[0']}\ot 1_{[1](1)}1_{[1']})\ot_A (1_A\ot 1_{[1](2)}).
\end{eqnarray*}
\equref{2.5.1} then follows from the left counit property in \equref{1.2.3}:
\begin{eqnarray*}
&&\hspace*{-2cm}
\rho(1_A)=(1_A\ot1_H)1_A=((\ul{\varepsilon}\circ \iota)\ot_A \ul{A\ot H})\ul{\Delta}((1_A\ot 1_H)1_A)\\
&=& (A\ot \epsilon)(1_{[0]} 1_{[0']}\ot 1_{[1](1)}1_{[1']})(1_A\ot 1_{[1](2)})\\
&=& \epsilon(1_{[1](1)}1_{[1']}) 1_{[0]} 1_{[0']}\ot 1_{[1](2)} \equal{\equref{2.2.1}} \epsilon(1_{[1]})1_{[0]}\ot 1_{[2]}.
\end{eqnarray*}
$\ul{2)\Rightarrow 1)}$. If conditions
(\ref{eq:2.1.1},\ref{eq:2.2.1},\ref{eq:2.2.2},\ref{eq:2.5.1}) are satisfied, then $\Delta$ is a coassociative
comultiplication on $A\ot H$. One equality in \equref{1.2.3} is equivalent to 
\equref{2.5.1}, and the other one can be proved as follows: we have shown in the proof of
\prref{2.3} that (\ref{eq:2.1.1}) and (\ref{eq:2.2.1}) imply that
$((A\ot H)\ot_A\ul{\varepsilon})\circ\Delta=\pi$, and
this entails that
$((\ul{A\ot H})\ot_A(\ul{\varepsilon}\circ \iota))\circ\ul{\Delta}=\ul{A\ot H}$.\\
$\ul{2)\Leftrightarrow 3)}$.
Using \equref{2.2.2}, we find that \equref{2.5.1} is equivalent to \equref{2.5.2}.\\
$\ul{2)\Leftrightarrow 4)}$. Using (\ref{eq:2.2.1}), we can prove the equivalence of \equref{2.5.1} and \equref{2.5.3}:
\begin{eqnarray*}
&&\hspace*{-1cm}
\epsilon(1_{[1]})1_{[0]}\ot 1_{[2]} \equal{\equref{2.2.1}}  \epsilon(1_{[1](1)}1_{[1']}) 1_{[0]} 1_{[0']}\ot 1_{[1](2)}= \epsilon(1_{[1']})1_{[0]} 1_{[0']}\ot 1_{[1]}.
\end{eqnarray*}
\end{proof}

\begin{proposition}\prlabel{2.6}
Let $A$ be a $k$-algebra, $H$ a $k$-bialgebra, and $\rho: A\to A\ot H$
a $k$-linear map. We call $(A,\rho)$ a (right) partial $H$-comodule algebra, and
say that $H$ coacts partially on $A$, if the following equivalent conditions are satisfied:\\
1) $(A\ot H,\Delta,\varepsilon)$ is a left unital lax $A$-coring;\\
2) $A$ is a (right) lax $H$-comodule algebra and
\begin{equation}\eqlabel{2.6.1}
\epsilon(1_{[1]})1_{[0]}=1_A;
\end{equation}
3) The conditions (\ref{eq:2.1.1},\ref{eq:2.2.1}) and \equref{2.2.3}
are satisfied, for all $a,b\in A$.
\end{proposition}

\begin{proof}
$\ul{1)\Rightarrow 2)}$. It follows from \leref{2.2} that \equref{2.2.3}
holds; taking $a=1$, we find \equref{2.6.1}. It is clear that $A$
is a lax $H$-comodule algebra.\\
$\ul{2)\Rightarrow 1)}$. From \equref{2.6.1} it follows that $\varepsilon=\ul{\varepsilon}$, and thus $(A\ot H,\Delta,\varepsilon=\ul{\varepsilon})$ is a left unital lax $A$-coring.\\
$\ul{2)\Rightarrow 3)}$. Combining \equref{2.6.1} and \equref{2.2.2}, we find
\equref{2.2.3}.\\
$\ul{3)\Rightarrow 2)}$. Taking $a=1$ in \equref{2.2.3}, we find \equref{2.6.1}.
\equref{2.2.2} then follows from  \equref{2.6.1} and  \equref{2.2.3}. Also \equref{2.5.2} (or \equref{2.5.3}) follows immediately from \equref{2.6.1}.
\end{proof}

\begin{proposition}\prlabel{2.7}
$A$ is an $H$-comodule algebra if and only if it is at the same time a partial and
weak $H$-comodule algebra. 
\end{proposition}

\begin{proof}
One implication is obvious. Conversely, if $A$ is a weak $H$-comodule algebra,
then, by \prref{2.3}, (\ref{eq:2.1.1},\ref{eq:2.3.1},\ref{eq:2.3.2}) are satisfied.
If $A$ is a partial $H$-comodule algebra, then \equref{2.2.3} holds, by
\prref{2.6}. Then \equref{2.1.2} can be shown as follows:
$\rho(1_A)\equal{\equref{2.3.1}}\epsilon(1_{[1]})1_{[0]}\ot 1_H
\equal{\equref{2.2.3}}1_A\ot 1_H$.
\end{proof}

\subsubsection*{The lax Koppinen smash product}
Let $A$ be a right weak (or lax) $H$-comodule algebra, i.e. $(\Cc=A\ot H,\Delta,\ul{\varepsilon})$ is a (left unital) weak (or lax) $A$-coring. Given the $k$-module isomorphism $${}^*\Cc={}_A\textrm{Hom}(A\ot H,A)\cong \textrm{Hom}(H,A), f\mapsto f\circ (\eta_A \ot H),$$ 
the (right unital) weak (or lax) $A$-ring structure on ${}^*\Cc$ (see \prref{1.5}) induces a (right unital) weak (or lax) $A$-ring structure on $\textrm{Hom}(H,A)$. It is given by the following formulas, for all $a,b\in A, h\in H$ and $f,g\in \textrm{Hom}(H,A)$:
$(afb)(h)=a_{[0]} f(ha_{[1]})b$,
$(f\#g)(h)=f(h_{(2)})_{[0]} g(h_{(1)} f(h_{(2)})_{[1]})$
and $\eta: A\rightarrow \textrm{Hom}(H,A)$, $\eta(a)(h)=\epsilon(ha_{[1]})a_{[0]}$.
$\textrm{Hom}(H,A)$ with this weak (or lax) $A$-ring structure 
will be called the weak (or lax) Koppinen smash product. It is usually denoted by $\#(H,A)$.
The left dual of the corresponding $A$-coring $\ul{\Cc}=\Cc1_A$ is then isomorphic to
$$\ul{\#}(H,A)=1_A\#(H,A)= \{f\in \#(H,A)~|~f(h)=1_{[0]} f(h1_{[1]}),~{\rm for~all}~h\in H\}.$$

\begin{remark}\relabel{2.8} 
All results in this Section can also be proved in the more general context of entwining structures (see \cite{BrzezinskiM}). These are triples $(A,C,\psi)$, where $A$ is a $k$-algebra, $C$ a $k$-coalgebra, and $\psi:C\ot A\to A\ot C$ a $k$-linear map satisfying some compatibility conditions. A right $H$-comodule algebra $A$ gives rise to an entwining structure $(A,H,\psi)$, where $\psi: H\ot A\to A\ot H$ is defined by $\psi(h\ot a)=a_{[0]}\ot ha_{[1]}$. The analogue of \prref{2.4} then recovers these entwining structures, whereas \prref{2.3} recovers the so-called weak entwining structures. These were first introduced by the first author and De Groot
\cite{CaenepeelDG00}. Wisbauer \cite{Wisbauer} introduced weak corings,
and obtained a one-to-one correspondence (see \cite[4.1]{Wisbauer}) between weak corings and weak entwining structures. In a remark following 4.1 in \cite{Wisbauer}, it is observed that the defining axioms of weak entwining structures in
\cite{CaenepeelDG00} and \cite{Wisbauer} are not the same. It follows from the analogue of \prref{2.3}
that the two sets of axioms are equivalent.
\end{remark}

\begin{example}\exlabel{3.2}
Let $e\in H$ be an idempotent such that $e\ot e=\Delta(e)(e\ot 1)$ and $\epsilon(e)=1$.
A partial $H$-coaction on $A=k$ is given by
$\rho(x)=x\ot e\in k\ot_k H$. (\ref{eq:2.1.1},\ref{eq:2.2.1},\ref{eq:2.2.3}) can be
verified easily:
$
\rho(x)\rho(y)= xy\ot e^2=xy\ot e=\rho(xy)$;
$\rho(x_{[0]})\ot x_{[1]}= x\ot e\ot e=x\ot e_{(1)}e\ot e_{(2)}=
x_{[0]}1_{[0]}\ot x_{[1](1)}1_{[1]}\ot x_{[1](2)}$;
$x\epsilon(e)=x$.\\
Such an idempotent $e$ exists in a finite dimensional semisimple Hopf algebra:
take a left integral $t$ such that $\epsilon(t)=1$. $t$ is an idempotent, since
$t^2=\epsilon(t)t=t$, and $\Delta(t)(t\ot 1)=t_{(1)}t\ot t_{(2)}=
\epsilon(t_{(1)})t\ot t_{(2)}=t\ot t$.
\end{example}

\begin{example}\exlabel{3.3}
Let $H$ be Sweedler's 4-dimensional Hopf algebra over a field $k$ with ${\rm char}(k)\neq 2$. Recall that, as a $k$-algebra, $H$ is generated by two elements $c$ and $x$, with relations $c^2=1$, $x^2=0$ and $xc=-cx$. Then $H$ is a 4-dimensional vector space with basis $\{1,c,x,cx\}$. The coalgebra structure is induced by $\Delta(c)=c\ot c$, $\Delta(x)=c\ot x+ x\ot 1$, $\epsilon(c)=1$ and $\epsilon(x)=0$. Note that $H$ is not semisimple.\\
We can define a partial coaction of $H$ on $A=k$ as in \exref{3.2}. For $\alpha\in k$,
consider $e=e_\alpha=\frac{1}{2}+\frac{1}{2}c+\alpha cx\in H$. Then it is easily checked that $e$ is an idempotent, such that $e\ot e=\Delta(e)(e\ot 1)$ and $\epsilon(e)=1$.\\
Now consider the subalgebra $B=k[x]$ of $H$. The map $\rho: B\to B\ot H$, given by
$\rho(1)=\frac{1}{2}\ot 1+\frac{1}{2}\ot c+\frac{1}{2}\ot cx$ and $ \rho(x)=\frac{1}{2}x\ot 1+\frac{1}{2}x\ot c+\frac{1}{2}x\ot cx$,
defines a partial coaction of $H$ on $B$.
\end{example}

\section{Lax relative Hopf modules}\selabel{2a}
Let $A$ be a right lax $H$-comodule algebra, $\Cc=A\ot H$ the associated lax $A$-coring,
and $\ul{\Cc}=(A\ot H)1_A$ the associated $A$-coring.
For a $k$-linear map
$\rho: M\to M\ot H$, we will adopt the notation
$\rho(m)=m_{[0]}\ot m_{[1]}$, $(\rho\ot H)(\rho(m))=\rho^2(m)=m_{[0]}\ot m_{[1]}\ot m_{[2]}$,
etc.

A lax relative Hopf module is a right $A$-module $M$, together with a $k$-linear map
$\rho: M\to M\ot H$ such that the following conditions are satisfied, for all $m\in M$:
\begin{eqnarray}
&&\hspace*{-2cm} m_{[0]}\epsilon(m_{[1]})=m;\eqlabel{2a.1.1}\\
&&\hspace*{-2cm} \rho^2(m)= m_{[0]}1_{[0]}\ot m_{[1](1)}1_{[1]}\ot m_{[1](2)}1_{[2]};
\eqlabel{2a.1.2}\\
&&\hspace*{-2cm}
\rho(ma)=m_{[0]}a_{[0]}\ot m_{[1]}a_{[1]}.\eqlabel{2a.1.3}
\end{eqnarray}
A morphism between two lax relative Hopf modules $M$ and $N$ is a right $A$-linear map
$f: M\to N$ such that $f(m_{[0]})\ot m_{[1]}=f(m)_{[0]}\ot f(m)_{[1]}$, for all
$m\in M$. $\Mm_A^H$ will denote the category of lax relative Hopf modules.

\begin{proposition}\prlabel{2a.1}
For a right lax $H$-comodule algebra $A$, the categories $\Mm^{\ul{\Cc}}$ and
$\Mm_A^H$ are isomorphic.
\end{proposition}

\begin{proof}
Let $M$ be a right $A$-module.
We will first show that $\Hom_A(M,M\ot_A\ul{\Cc})$ is isomorphic to the submodule of
$\Hom(M,M\ot H)$ consisting of maps $\rho$ satisfying \equref{2a.1.3}.
For $\rho: M\to M\ot H$ satisfying \equref{2a.1.3}, consider the map $\alpha(\rho):
M\to M\ot_A \ul{\Cc}$,
$\alpha(\rho)(m)=m_{[0]}\ot_A (1_{[0]}\ot m_{[1]}1_{[1]})$.
$\alpha(\rho)$ is right $A$-linear since
\begin{eqnarray*}
&&\hspace*{-1cm}
\alpha(\rho)(ma)=m_{[0]}a_{[0]}\ot_A (1_{[0]}\ot m_{[1]}a_{[1]}1_{[1]})=
 m_{[0]}\ot_A (a_{[0]}1_{[0]}\ot m_{[1]}a_{[1]}1_{[1]})\\
&\equal{\equref{2.1.1}}& m_{[0]}\ot_A (a_{[0]}\ot m_{[1]}a_{[1]})
= m_{[0]}\ot_A (1_{[0]}\ot m_{[1]}1_{[1]})a=(\alpha(\rho)(m))a.
\end{eqnarray*}
Conversely, take $\tilde{\rho}\in \Hom_A(M,M\ot_A\ul{\Cc})$, and define $\beta(\tilde{\rho})
\in \Hom(M, M\ot H)$ as follows: for $m\in M$, there exist  (a finite number of) $m_i\in M$
and $h_i\in H$ such that $\tilde{\rho}(m)= \sum_i m_i\ot_A (1_{[0]}\ot h_i1_{[1]})$; then we
define 
$\beta(\tilde{\rho})(m)=\sum_i m_i1_{[0]}\ot h_i1_{[1]}$. We claim that
$\beta(\tilde{\rho})$ satisfies \equref{2a.1.3}. Indeed,
$\tilde{\rho}(ma)= \sum_i m_i\ot_A (a_{[0]}\ot h_ia_{[1]})
=\sum_i m_ia_{[0]}\ot_A (1_{[0]}\ot h_ia_{[1]}1_{[1]})$, hence
$\beta(\tilde{\rho})(ma)=\sum_i m_ia_{[0]}1_{[0]}\ot h_ia_{[1]}1_{[1]}=
\sum_i m_i1_{[0]}a_{[0]}\ot h_i1_{[1]}a_{[1]}$.
$\alpha$ and $\beta$ are inverses since
$\beta(\alpha(\rho))(m)= m_{[0]}1_{[0]} \ot m_{[1]}1_{[1]}=\rho(m)$ and
$
\alpha(\beta(\tilde{\rho}))(m)=\sum_i m_i1_{[0]} \ot_A (1_{[0']}\ot h_i1_{[1]}1_{[1']})
= \sum_i m_i\ot_A (1_{[0]}\ot h_i1_{[1]})=\tilde{\rho}(m)$.\\
Now take $\rho: M\to M\ot H$ satisfying \equref{2a.1.3} and the corresponding
right $A$-linear map $\tilde{\rho}$. We claim that $\tilde{\rho}$ is coassociative if and only
if $\rho$ satisfies \equref{2a.1.2}. First compute
$$\tilde{\rho}^2(m)=
m_{[0]}\ot_A (1_{[0']}\ot m_{[1]}1_{[1']})\ot_A (1_{[0]}\ot m_{[2]}1_{[1]});$$
$$((M\ot_A\ul{\Delta})\circ \tilde{\rho})(m)=
m_{[0]}\ot_A (1_{[0]}\ot m_{[1](1)}1_{[1](1)})\ot_A (1_A\ot m_{[1](2)}1_{[1](2)}).$$
If $\tilde{\rho}$ is coassociative, then it follows that
\begin{eqnarray*}
&&\hspace*{-2cm}
m_{[0]}1_{[0']}1_{[0]}\ot m_{[1]}1_{[1']}1_{[1]}\ot m_{[2]}1_{[2]}
\equal{\equref{2.1.1}}
m_{[0]}1_{[0]}\ot m_{[1]}1_{[1]}\ot m_{[2]}1_{[2]}\\
&\equal{\equref{2a.1.3}}& \rho(m_{[0]}1_{[0]})\ot m_{[1]}1_{[1]}\equal{\equref{2a.1.3}}\rho(m_{[0]})\ot m_{[1]}
\end{eqnarray*}
equals
$$
m_{[0]}1_{[0]}1_{[0']}\ot m_{[1](1)}1_{[1](1)}1_{[1']}\ot m_{[1](2)}1_{[1](2)}
\equal{\equref{2.2.1}} m_{[0]}1_{[0]}\ot m_{[1](1)}1_{[1]} \ot m_{[1](2)}1_{[2]} ,
$$
and  \equref{2a.1.2} follows. Conversely, if \equref{2a.1.2} holds, then
\begin{eqnarray*}
&&\hspace*{-1.5cm}
((M\ot_A\ul{\Delta})\circ \tilde{\rho})(m)\\
&=&
m_{[0]}\ot_A (1_{[0]} 1_{[0']} \ot m_{[1](1)}1_{[1](1)}1_{[1']})
\ot_A (1_A\ot m_{[1](2)}1_{[1](2)})\cdot 1_A\\
&\equal{\equref{2.2.1}}&
m_{[0]}\ot_A (1_{[0]}  \ot m_{[1](1)}1_{[1]})
\ot_A (1_A\ot m_{[1](2)}1_{[2]})\cdot 1_A\\
&=&
m_{[0]}1_{[0]} \ot_A (1_A \ot m_{[1](1)}1_{[1]})\cdot 1_A
\ot_A (1_A\ot m_{[1](2)}1_{[2]})\cdot 1_A\\
&\equal{\equref{2a.1.2}}&
m_{[0]}\ot_A (1_A\ot m_{[1]})\cdot 1_A\ot_A (1_A\ot m_{[2]})\cdot 1_A=\tilde{\rho}^2(m),
\end{eqnarray*}
so $\tilde{\rho}$ is coassociative. Finally, $\tilde{\rho}$ satisfies the counit property if
and only if
$m=m_{[0]}\epsilon(m_{[1]}1_{[1]})1_{[0]} = m_{[0]}\epsilon(m_{[1]})$.
\end{proof}

\section{Partial module algebras}\selabel{4}
Let $A$ be a $k$-algebra, $H$ a $k$-bialgebra, and $\kappa:H \ot A\to A,~\kappa(h\ot a)=h\cdot a$ a $k$-linear map.
We assume that $A\ot H$ is a left unital $A$-bimodule under the action
$b'(a\ot h)b=b'a(h_{(1)}\cdot b)\ot h_{(2)}$. Then the following condition is satisfied, for all $a,b\in A$
and $h\in H$ (see \equref{2.1.1}):
\begin{equation}\eqlabel{4.1.0}
h\cdot(ab)=(h_{(1)}\cdot a)(h_{(2)}\cdot b).
\end{equation}
The map $\mu: (A\ot H)\ot_A(A\ot H)\to A\ot H$,
$\mu((a\ot h)\ot_A (b\ot g))=a(h_{(1)}\cdot b)\ot h_{(2)}g$,
is well-defined since
\begin{eqnarray*}
&&\hspace{-2cm}
\mu((a\ot h)a'\ot_A (b\ot g))=\mu((a(h_{(1)}\cdot a')\ot h_{(2)})\ot_A (b\ot g))\\
&=&a(h_{(1)}\cdot a')(h_{(2)}\cdot b)\ot h_{(3)}g
\equal{\equref{4.1.0}}a(h_{(1)}\cdot (a'b))\ot h_{(2)}g\\
&=&\mu((a\ot h)\ot_A (a'b\ot g)).
\end{eqnarray*}
$A\# H$ will be our notation for $A\ot H$ together with the multiplication $\mu$.
We then write
$\mu((a\ot h)\ot_A (b\ot g))=(a\# h)(b\#g)$.
We also consider the maps
\begin{eqnarray*}
\eta&=&A\ot \eta_H: A\to A\ot H,~~\eta(a)=a\ot 1_H;\\
\ul{\eta}&=&\pi\circ \eta: A\to \ul{A\ot H},~~\ul{\eta}(a)=a(1_H\cdot 1_A)\ot 1_H.
\end{eqnarray*}

\begin{lemma}\lelabel{4.1}
Assume that $\kappa:H\ot A\to A$ satisfies \equref{4.1.0}.\\
1)
$\mu$ is right $A$-linear if and only if $\mu$ is associative if and only if
\begin{equation}\eqlabel{4.1.1}
h\cdot (a(g\cdot b))=(h_{(1)}\cdot a)((h_{(2)}g)\cdot b),
\end{equation}
for all $a,b\in A$ and $h,g\in H$.\\
2)
$\eta$ is right $A$-linear if and only if, for all $a\in A$,
\begin{equation}\eqlabel{4.1.2}
1_H\cdot a=a.
\end{equation}
3)
$\ul{\eta}$ is right $A$-linear if and only if, for all $a\in A$,
\begin{equation}\eqlabel{4.1.3}
a(1_H\cdot 1_A)=1_H\cdot a.
\end{equation}
\end{lemma}

\begin{proof}
1) $\mu$ is right $A$-linear if and only if
$$(a'\# h)((a\# g)b)=(a'\#h)(a(g_{(1)}\cdot b) \#g_{(2)})=a'(h_{(1)}\cdot (a(g_{(1)}\cdot b)))  \# h_{(2)}g_{(2)}$$
equals
$$((a'\# h)(a\# g))b=(a'(h_{(1)}\cdot a)\#h_{(2)}g)b=a'(h_{(1)}\cdot a)((h_{(2)}g_{(1)})\cdot b)\# h_{(3)}g_{(2)},$$
for all $a',a,b\in A$ and $h,g\in H$. This is equivalent to \equref{4.1.1}. It is obvious that $\mu$ is associative if and only if $\mu$ is right $A$-linear.\\
2) $\eta$ is right $A$-linear if and only if $\eta(a)=a\ot 1_H$ equals
$\eta(1_A)a=(1_A\ot 1_H)a=1_H\cdot a \ot 1_H$, which is equivalent to \equref{4.1.2}.\\
3) $\ul{\eta}$ is right $A$-linear if and only if $\ul{\eta}(a)=a(1_H\cdot 1_A)\ot 1_H$ equals
$$\ul{\eta}(1_A)a=(1_H\cdot 1_A\ot 1_H)a=(1_H\cdot 1_A)(1_H\cdot a)\ot 1_H\equal{\equref{4.1.0}}1_H\cdot a\ot 1_H,$$
and this is equivalent to \equref{4.1.3}.
\end{proof}

\begin{proposition}\prlabel{4.2}
Let $A$ be a $k$-algebra, $H$ a $k$-bialgebra and $\kappa: H\ot A\to A$ a $k$-linear map.
$A$ is termed a (left) weak $H$-module algebra if the following equivalent conditions
are satisfied:\\
1) $(A\# H,\mu,\ul{\eta})$ is a left unital weak $A$-ring;\\
2) the conditions (\ref{eq:4.1.0},\ref{eq:4.1.1},\ref{eq:4.1.3},\ref{eq:4.2.1})
are satisfied, for all $a,b\in A$ and $h,g\in H$;\\
3) the conditions 
(\ref{eq:4.1.0},\ref{eq:4.1.3},\ref{eq:4.2.1},\ref{eq:4.2.2}) are satisfied, for all $a,b\in A$ and $h,g\in H$.
\begin{eqnarray}
h\cdot 1_A&=&\epsilon(h)1_H\cdot 1_A;\eqlabel{4.2.1}\\
h\cdot (g\cdot a)&=&(hg)\cdot a.\eqlabel{4.2.2}
\end{eqnarray}
\end{proposition}

\begin{proof}
$\ul{1)\Rightarrow 2)}$. (\ref{eq:4.1.0},\ref{eq:4.1.1},\ref{eq:4.1.3})
follow from \leref{4.1}. From \equref{1.4.2}, it follows that
$(1_A\# h)1_A= (h_{(1)}\cdot 1_A)\#h_{(2)}$
equals
$
((1_H\cdot 1_A)\#1_H)(1_A\# h)=(1_H\cdot 1_A)(1_H\cdot 1_A)\# h\equal{\equref{4.1.0}}(1_H\cdot 1_A)\#h,$
and \equref{4.2.1} follows by applying $A\#\epsilon$.\\
$\ul{2)\Rightarrow 1)}$. If (\ref{eq:4.1.0},\ref{eq:4.1.1},\ref{eq:4.1.3}) hold, then
we only need to verify \equref{1.4.2}, by \leref{4.1}. We compute that
\begin{eqnarray*}
&&\hspace*{-1.2cm}
((1_H\cdot 1_A)\# 1_H)(a\# h)=(1_H\cdot 1_A)(1_H\cdot a)\#h\equal{\equref{4.1.0}}
(1_H\cdot a)\# h
\equal{\equref{4.1.3}} a(1_H\cdot 1_A)\#h\\
&=&a\epsilon(h_{(1)})(1_H\cdot 1_A)\#h_{(2)}\equal{\equref{4.2.1}}
a(h_{(1)}\cdot 1_A)\#h_{(2)}=(a\# h)1_A;\\
&&\hspace*{-1cm}
(a\# h)((1_H\cdot 1_A)\# 1_H)=a(h_{(1)}\cdot (1_H\cdot 1_A))\# h_{(2)}\\
&\equal{\equref{4.1.1}}&a(h_{(1)}\cdot 1_A)(h_{(2)}\cdot 1_A)\#h_{(3)}
\equal{\equref{4.1.0}}a(h_{(1)}\cdot 1_A)\# h_{(2)}=(a\# h)1_A.
\end{eqnarray*}
$\ul{2)\Rightarrow 3)}$. We have that
\begin{eqnarray*}
&&\hspace*{-2cm}
h\cdot (g\cdot a)\equal{\equref{4.1.1}} (h_{(1)}\cdot 1_A)((h_{(2)}g)\cdot a)
\equal{\equref{4.2.1}} \epsilon(h_{(1)})(1_H\cdot 1_A)((h_{(2)}g)\cdot a)\\
&=& (1_H\cdot 1_A)((hg)\cdot a)=\epsilon(h_{(1)}g_{(1)})(1_H\cdot 1_A)((h_{(2)}g_{(2)})\cdot a)\\
&\equal{\equref{4.2.1}}&((h_{(1)}g_{(1)})\cdot 1_A)((h_{(2)}g_{(2)})\cdot a)\equal{\equref{4.1.0}}(hg)\cdot a  ,
\end{eqnarray*}
and \equref{4.2.2} follows.\\
$\ul{3)\Rightarrow 2)}$. We have to show that \equref{4.1.1} holds. Indeed,
$
h\cdot(a(g\cdot b))\equal{\equref{4.1.0}}(h_{(1)}\cdot a)(h_{(2)}\cdot (g\cdot b))
\equal{\equref{4.2.2}} (h_{(1)}\cdot a)((h_{(2)}g)\cdot b)$.
\end{proof}

\begin{proposition}\prlabel{4.3}
Let $A$ be a $k$-algebra, $H$ a $k$-bialgebra and $\kappa: H\ot A\to A$ a $k$-linear map.
The following assertions are equivalent:

1) $(A\# H,\mu,{\eta})$ is an $A$-ring;\\
2) $(A\# H,\mu,{\eta})$ is a left unital weak $A$-ring;\\
3) $A$ is a left $H$-module algebra; this means that the following conditions hold, for all $a,b\in A$ and $h,g\in H$: (\ref{eq:4.1.0},\ref{eq:4.1.2},\ref{eq:4.2.2}) and
\begin{equation}\eqlabel{4.3.1}
h\cdot 1_A=\epsilon(h)1_A.
\end{equation}
\end{proposition}

\begin{proof}
$\ul{1)\Rightarrow 2)}$ is trivial; $\ul{3)\Rightarrow 1)}$ is well-known (and easy to prove).\\
$\ul{2)\Rightarrow 3)}.$ It follows from \leref{4.1} that (\ref{eq:4.1.0},\ref{eq:4.1.1},\ref{eq:4.1.2})
hold. Using \equref{1.4.2}, we find that
$
h_{(1)}\cdot 1_A\# h_{(2)}=(1_A\# h)1_A=((1_H\cdot 1_A)\#1_H)(1_A\#h)
=(1_H\cdot 1_A)(1_H\cdot 1_A)\#h \equal{\equref{4.1.0}}(1_H\cdot 1_A)\#h \equal{\equref{4.1.2}}1_A\# h$,
hence \equref{4.3.1} follows after we apply $A\# \epsilon$ to both sides. Taking $a=1$ in \equref{4.1.1} and using \equref{4.3.1},
we find \equref{4.2.2}.
\end{proof}

\begin{proposition}\prlabel{4.4}
Let $A$ be a $k$-algebra, $H$ a $k$-bialgebra and $\kappa: H\ot A\to A$ a $k$-linear map.
$A$ is called a (left) lax $H$-module algebra if the following equivalent properties are
satisfied:\\
1) $(A\# H,\mu,\ul{\eta})$ is a left unital lax $A$-ring;\\
2) the conditions (\ref{eq:4.1.0},\ref{eq:4.1.1},\ref{eq:4.1.3},\ref{eq:4.4.1}) 
are fulfilled, for all $a,b\in A$ and $h,g\in H$;\\
3) the conditions (\ref{eq:4.1.0},\ref{eq:4.1.1},\ref{eq:4.1.3},\ref{eq:4.4.2}) 
are fulfilled, for all $a,b\in A$ and $h,g\in H$.
\begin{eqnarray}
a(h\cdot 1_A)&=&1_H\cdot (a(h\cdot 1_A))\eqlabel{4.4.1};\\
a(h\cdot 1_A)&=&(1_H\cdot a)(h\cdot 1_A)\eqlabel{4.4.2}.
\end{eqnarray}
\end{proposition}

\begin{proof}
$\ul{1)\Rightarrow 2)}$. (\ref{eq:4.1.0},\ref{eq:4.1.1},\ref{eq:4.1.3})
follow from \leref{4.1}. Using \equref{1.4.1}, we find
\begin{eqnarray*}
&&\hspace*{-1.5cm}
a(h_{(1)}\cdot 1_A)\# h_{(2)}= (a\#h)1_A=((1_H\cdot 1_A)\#1_H)((a\#h)1_A)\\
&=& ((1_H\cdot 1_A)\#1_H)(a(h_{(1)}\cdot 1_A)\# h_{(2)})\\
&=& (1_H\cdot 1_A)(1_H\cdot (a(h_{(1)}\cdot 1_A)))\#h_{(2)} \equal{\equref{4.1.0}}(1_H\cdot (a(h_{(1)}\cdot 1_A)))\#h_{(2)},
\end{eqnarray*}
so \equref{4.4.1} follows after applying $A\#\epsilon$.\\
$\ul{2)\Rightarrow 3)}$ follows immediately by \equref{4.1.1}.\\
$\ul{3)\Rightarrow 1)}$. It follows from the above computations that
\begin{eqnarray*}
&&\hspace*{-1.5cm}
((1_H\cdot 1_A)\#1_H)((a\#h)1_A)= (1_H\cdot (a(h_{(1)}\cdot 1_A)))\#h_{(2)}\\
&\equal{\equref{4.1.1}} &(1_H\cdot a)(h_{(1)}\cdot 1_A)\#h_{(2)} \equal{\equref{4.4.2}} a(h_{(1)}\cdot 1_A)\# h_{(2)}= (a\#h)1_A.
\end{eqnarray*} 
We also have
\begin{eqnarray*}
&&\hspace*{-2cm}
((a\#h)1_A)((1_H\cdot 1_A)\#1_H)=(a(h_{(1)}\cdot 1_A)\# h_{(2)})((1_H\cdot 1_A)\#1_H)\\
&=& a(h_{(1)}\cdot 1_A)(h_{(2)}\cdot (1_H\cdot 1_A))\#h_{(3)}\\
&\equal{\equref{4.1.1}}& a(h_{(1)}\cdot 1_A)(h_{(2)}\cdot 1_A)(h_{(3)}\cdot 1_A)\#h_{(4)}\\
&\equal{\equref{4.1.0}}& a(h_{(1)}\cdot 1_A)\#h_{(2)}
= (a\#h)1_A,
\end{eqnarray*}
as needed.
\end{proof}

\begin{proposition}\prlabel{4.5}
Let $A$ be a $k$-algebra, $H$ a $k$-bialgebra and $\kappa: H\ot A\to A$ a $k$-linear map.
We say that $A$ is a (left) partial $H$-module algebra, or that $H$ acts partially on $A$,
if the following assertions are satisfied:\\
1) $(A\# H,\mu,{\eta})$ is a left unital lax $A$-ring;\\
2) $A$ is a left lax $H$-module algebra and
$1_H\cdot 1_A=1_A$;\\
3) the conditions (\ref{eq:4.1.0},\ref{eq:4.1.1},\ref{eq:4.1.2}) are satisfied.
\end{proposition}

\begin{proof}
$\ul{1)\Rightarrow 2)}$. Let $(A\# H,\mu,{\eta})$ be a left unital lax $A$-ring.
It follows from \leref{4.1} that \equref{4.1.2} holds. Taking
$a=1$ in \equref{4.1.2}, we find that $1_H\cdot 1_A=1_A$. This implies that $\eta=\ul{\eta}$, so
$(A\# H,\mu,\ul{\eta})$ is a left unital lax $A$-ring.\\
$\ul{2)\Rightarrow 1)}$ follows also from the fact that $1_H\cdot 1_A=1_A$ implies that $\eta=\ul{\eta}$.\\
$\ul{1)\Rightarrow 3)}$ follows immediately from \leref{4.1}.\\
$\ul{3)\Rightarrow 1)}$. We have to show that \equref{1.4.1} holds:
\begin{eqnarray*}
&&\hspace*{-1.5cm}
((a\#h)1_A)((1_H\cdot 1_A)\#1_H)= a(h_{(1)}\cdot 1_A)(h_{(2)}\cdot (1_H\cdot 1_A))\#h_{(3)}\\
&\equal{\equref{4.1.2}}& a(h_{(1)}\cdot 1_A)(h_{(2)}\cdot 1_A)\#h_{(3)}\equal{\equref{4.1.0}} a(h_{(1)}\cdot 1_A)\#h_{(2)}=(a\#h)1_A;\\
&&\hspace*{-1.5cm}
((1_H\cdot 1_A)\#1_H)((a\#h)1_A)= (1_H\cdot (a(h_{(1)}\cdot 1_A)))\#h_{(2)}\\
&\equal{\equref{4.1.2}}& a(h_{(1)}\cdot 1_A)\#h_{(2)}=(a\#h)1_A.
\end{eqnarray*}
\end{proof}

\begin{proposition}
$A$ is an $H$-module algebra if and only if it is at the same time a weak and partial $H$-module algebra.
\end{proposition}

\begin{proof}
One implication is obvious. Conversely, if $A$ is a weak $H$-module algebra, then (\ref{eq:4.1.0},\ref{eq:4.1.3},\ref{eq:4.2.1},\ref{eq:4.2.2}) are satisfied (cf. \prref{4.2}). If $A$ is a partial $H$-module algebra, then (\ref{eq:4.1.2}) holds (cf. \prref{4.5}) and (\ref{eq:4.3.1}) can be verified 
using \equref{4.2.1} and $1_H\cdot 1_A=1_A$:
$h\cdot 1_A= \epsilon(h)1_H\cdot 1_A  =  \epsilon(h)1_A$.
\end{proof}

\begin{theorem}\thlabel{4.8}
Let $A$ be a $k$-algebra, and $H$ a finitely generated projective $k$-bialgebra. Then $A$ is a (right) lax (resp. partial) $H$-comodule algebra if and only if $A^{\rm op}$ is a (left) lax (resp. partial) $H^{*{\rm cop}}$-module algebra.
\end{theorem}

\begin{proof}
Let $\{h_i,h_i^* ~ \vert ~ i=1,\ldots,n\}$ be a dual basis for $H$. Then it is well-known
(see for example \cite[(1.5)]{CMZ}) that
\begin{equation}\eqlabel{4.6.1}
\sum_i \Delta(h_i)\ot h_i^* = \sum_{i,j}h_i\ot h_j \ot h_i^* * h_j^*.
\end{equation}
We also have that 
$\Hom(A, A\ot H)\cong \Hom(H^{*}\ot A, A)$ as $k$-modules.
The isomorphism can be described as follows. If $\kappa:H^*\ot A\to A,\kappa(h^*\ot a)=h^*\leftact a$
corresponds to $\rho:A\to A\ot H,\rho(a)=a_{[0]}\ot a_{[1]}$, then
\begin{equation}\eqlabel{4.6.2}
h^*\leftact a=h^*(a_{[1]})a_{[0]}~~{\rm and}~~
\rho(a)=\sum_{i} h_i^*\leftact a \ot h_i.
\end{equation}
Assume that $A$ is a right lax $H$-comodule algebra; we will show that
$A^{\rm op}$ is a left lax $H^{*{\rm cop}}$-module algebra. The multiplication in $A^{\rm op}$
will be denoted by a dot: $a\cdot b=ba$. We have to show that $\kappa$ satisfies
(\ref{eq:4.1.0},\ref{eq:4.1.1},\ref{eq:4.1.3}) and \equref{4.4.1}.
\begin{eqnarray*}
&&\hspace*{-1.2cm}
h^*\leftact (a\cdot b)= h^*((ba)_{[1]})(ba)_{[0]}\equal{\equref{2.1.1}}h^*(b_{[1]}a_{[1]})b_{[0]}a_{[0]}\\
&=& h^*_{(1)}(b_{[1]})h^*_{(2)}(a_{[1]})b_{[0]}a_{[0]}= (h^*_{(2)}\leftact a)\cdot (h^*_{(1)}\leftact b);\\
&&\hspace*{-1.2cm}
h^*\leftact (a\cdot (g^*\leftact b))= g^*(b_{[1]})h^*\leftact (b_{[0]}a)= g^*(b_{[2]})h^*(b_{[1]}a_{[1]})b_{[0]}a_{[0]}\\
&\equal{\equref{2.2.1}} &  g^*(b_{[1](2)})h^*(b_{[1](1)}1_{[1]}a_{[1]})b_{[0]}1_{[0]}a_{[0]}\equal{\equref{2.1.1}}g^*(b_{[1](2)})h^*(b_{[1](1)}a_{[1]})b_{[0]}a_{[0]}\\
&=& g^*(b_{[1](2)})h^*_{(1)}(b_{[1](1)})h^*_{(2)}(a_{[1]})b_{[0]}a_{[0]}\\
&=&(h^*_{(1)}*g^*)(b_{[1]})h^*_{(2)}(a_{[1]})b_{[0]}a_{[0]}\\
&=& (h^*_{(2)}\leftact a)\cdot ((h^*_{(1)}*g^*)\leftact b);\\
&&\hspace*{-1.2cm}
a\cdot (\epsilon \leftact 1_A)=\epsilon(1_{[1]})1_{[0]}a
\equal{\equref{2.2.2}} \epsilon(a_{[1]})a_{[0]}=\epsilon\leftact a;\\
&&\hspace*{-1.2cm}
\epsilon\leftact (a\cdot (h^*\leftact 1_A))= h^*(1_{[1]})\epsilon \leftact (1_{[0]}a)= h^*(1_{[2]})\epsilon(1_{[1]}a_{[1]})1_{[0]}a_{[0]}\\
&=&h^*(1_{[2]})\epsilon(1_{[1]})\epsilon(a_{[1]})1_{[0]}a_{[0]}\equal{\equref{2.2.2}} h^*(1_{[1]})\epsilon(1_{[1']})\epsilon(a_{[1]})1_{[0']}1_{[0]}a_{[0]}\\
&\equal{\equref{2.5.2}}& h^*(1_{[1]})\epsilon(a_{[1]})1_{[0]}a_{[0]}\equal{\equref{2.2.2}} h^*(1_{[1]})\epsilon(1_{[1']})1_{[0]}1_{[0']}a\\
&\equal{\equref{2.5.3}}& h^*(1_{[1]})1_{[0]}a=a\cdot (h^*\leftact 1_A).
\end{eqnarray*}
If $A$ is a right partial $H$-comodule algebra, then $A^{\rm op}$ is left partial $H^{*{\rm cop}}$-module algebra. It suffices to observe that, by \equref{2.6.1},
$\epsilon\leftact 1_A= \epsilon(1_{[1]})1_{[0]} = 1_A$.\\
Conversely, if $A^{\rm op}$ is a lax, resp. partial $H^{*{\rm cop}}$-module algebra,
then $A$ is a lax, resp. partial $H$-comodule algebra. The computations are similar
to the ones above, and are left to the reader.
\end{proof}

\begin{proposition}\prlabel{4.10}
Assume that $H$ is finitely generated and projective as a $k$-module.
Let $A$  be a right lax $H$-comodule algebra, and consider the left lax $H^{*{\rm cop}}$-module algebra structure on $A^{\rm op}$.
Then ${}^*(A\ot H)^{\rm op}$ is isomorphic to $A^{\rm op}\# H^{*{\rm cop}}$ as left unital lax $A^{\rm op}$-rings,
and ${}^*(\ul{A\ot H})^{\rm op}$ is isomorphic to $\ul{A^{\rm op}\# H^{*{\rm cop}}}$ as
$A^{\rm op}$-rings. Consequently the categories ${}_{\ul{A^{\rm op}\# H^{*\rm{cop}}}}\Mm$ and $\Mm_A^H$
are isomorphic.
\end{proposition}

\begin{proof}
We know from \seref{2} that ${}^*(A\ot H)^{\rm op}\cong \#(H,A)^{\rm op}$, with multiplication
$$(f\bullet g)(h)=g(h_{(2)})_{[0]} f(h_{(1)}g(h_{(2)})_{[1]})$$
in $\#(A,H)^{\rm op}$. The map
$$\alpha: A^{\rm op}\ot H^{*{\rm cop}}\to \Hom(H,A),~~\alpha(a\# h^*)(h)=a h^*(h)$$
is an isomorphism of $k$-modules since $H$ is finitely generated and projective. The first
statement follows after we show that $\alpha$ preserves the multiplication.
\begin{eqnarray*}
&&\hspace*{-2cm}
\alpha((a\# h^*)(b\#g^*))(h)=
\alpha ( a\cdot (h^*_{(2)}\leftact b)\# h^*_{(1)}*g^*)(h)\\
&=&(h^*_{(2)}\leftact b)a  (h^*_{(1)}*g^*)(h)= h^*_{(2)}(b_{[1]})b_{[0]}a h^*_{(1)}(h_{(1)})g^*(h_{(2)})\\
&=& h^*(h_{(1)}b_{[1]})g^*(h_{(2)})b_{[0]}a= g^*(h_{(2)})b_{[0]}\alpha(a\#h^*)(h_{(1)}b_{[1]})\\
&=& (\alpha(b\#g^*)(h_{(2)}))_{[0]} \alpha(a\# h^*)(h_{(1)}(\alpha(b\#g^*)(h_{(2)}))_{[1]})\\
&=&(\alpha(a\# h^*)\bullet \alpha(b\#g^*))(h).
\end{eqnarray*}
Applying \prref{1.5} we see that
${}^*(\ul{A\ot H})^{\rm op}\cong (1_A {}^*(A\ot H))^{\rm op}=$\\
${}^*(A\ot H)^{\rm op} 1_A\cong (A^{\rm op}\# H^{*{\rm cop}})1_A=
\ul{A^{\rm op}\# H^{*{\rm cop}}}$.
\end{proof}

We will investigate the notion of partial $H$-action in the particular situation where $H=kG$ is a group algebra, and the one where $H=U(L)$ is the universal enveloping algebra of a Lie algebra.
In the first case we will recover the partial group actions introduced in \cite{DFP}, at least in the
case where the involved ideals are generated by idempotents. In this particular situation, 
studied in \cite{DFP}, partial Galois theory can be developed.

\subsubsection*{Partial group actions}
Let $G$ be a group, and $A$ a $k$-algebra. A {\sl partial action} of $G$ on $A$
consists of a set of idempotents $\{e_\sigma~|~\sigma\in G\}\subset A$, and a set
of isomorphisms $\alpha_\sigma: e_{\sigma^{-1}}A\to e_\sigma A$ such that
$e_1=1_A$, $\alpha_1=A$
and
\begin{eqnarray}
e_\sigma\alpha_{\sigma\tau}(e_{\tau^{-1}\sigma^{-1}}a)&=&
\alpha_\sigma(e_{\sigma^{-1}}\alpha_\tau(e_{\tau^{-1}}a))\eqlabel{5.2.1};\\
\alpha_\sigma(e_{\sigma^{-1}}ab)&=&\alpha_\sigma(e_{\sigma^{-1}}a)\alpha_\sigma(e_{\sigma^{-1}}b);\eqlabel{5.2.2}\\
\alpha_\sigma(e_{\sigma^{-1}}) &=& e_\sigma,\eqlabel{5.2.2b}
\end{eqnarray}
for all $\sigma,\tau\in G$ and $a,b\in A$. This slightly generalizes the definitions in
\cite{DFP} and  \cite{CaenepeelDG05}: in \cite{DFP}, it is assumed that $A$ is commutative and that the isomorphisms $\alpha_\sigma$ are multiplicative;
in \cite{CaenepeelDG05}, it is assumed that, for all $\sigma\in G$, $e_\sigma$ is central and $\alpha_\sigma$ is multiplicative.
In both cases \equref{5.2.2} and \equref{5.2.2b} are automatically satisfied.

\begin{proposition}\prlabel{5.2}
Let $A$ be a $k$-algebra, and $G$ a group. Then there is a bijective correspondence
between partial $G$-actions and partial $kG$-actions on $A$.
\end{proposition}

\begin{proof}
Assume first that $kG$ acts partially on $A$. For each $\sigma\in G$, let
$e_\sigma=\sigma\cdot 1_A$.
Taking $a=c=1_A$ in \equref{4.1.0} we find that $e_\sigma^2=e_\sigma$. It follows
from \equref{4.1.2} that $e_1=1\cdot 1_A=1_A$. From \equref{4.1.1}, it follows that
\begin{equation}\eqlabel{5.2.3}
\sigma\cdot (\tau\cdot a)=e_\sigma ((\sigma\tau)\cdot a).
\end{equation}
We then compute
\begin{equation}\eqlabel{5.2.3b}
\sigma\cdot e_{\sigma^{-1}}=\sigma\cdot (\sigma^{-1}\cdot 1_A)=e_\sigma 1_A
=e_\sigma,
\end{equation}
and
\begin{equation}\eqlabel{5.2.4}
\sigma\cdot (e_{\sigma^{-1}}a)\equal{\equref{4.1.0}}
(\sigma\cdot e_{\sigma^{-1}})(\sigma\cdot a)=e_\sigma (\sigma\cdot a).
\end{equation}
It follows that the map $A\to A$, $a\mapsto \sigma\cdot a$ restricts to a map
$\alpha_\sigma: e_{\sigma^{-1}}A\to e_\sigma A$.
Observe that
\begin{equation}\eqlabel{5.2.5}
\sigma\cdot(e_{\sigma^{-1}}a)\equal{\equref{5.2.4}}e_\sigma(\sigma\cdot a)=
(\sigma\cdot 1_A)(\sigma\cdot a)\equal{\equref{4.1.0}}\sigma\cdot a.
\end{equation}
Now
$$
\alpha_{\sigma^{-1}}(\alpha_{\sigma}(e_{\sigma^{-1}}a)
\equal{\equref{5.2.4}}\alpha_{\sigma^{-1}}(e_\sigma (\sigma\cdot a))
\equal{\equref{5.2.5}}\sigma^{-1}\cdot(\sigma\cdot a)\\
\equal{\equref{5.2.3}}e_{\sigma^{-1}}((\sigma^{-1}\sigma)\cdot a)
\equal{\equref{4.1.2}}e_{\sigma^{-1}}a.
$$
In a similar way, we find that $\alpha_\sigma(\alpha_{\sigma^{-1}}(e_\sigma a))=
e_\sigma a$, and it follows that $\alpha_\sigma: e_{\sigma^{-1}}A\to e_\sigma A$
is an isomorphism. It is also clear that
$\alpha_1(a)=1\cdot a= a$ (use \equref{4.1.2})
and 
$
\alpha_\sigma(e_{\sigma^{-1}})=\sigma \cdot (e_{\sigma^{-1}})=e_\sigma
$ (use \equref{5.2.3b}).
\equref{5.2.2} can be shown as follows:
\begin{eqnarray*}
&&\hspace*{-2cm}
\alpha_\sigma(e_{\sigma^{-1}}a)\alpha_\sigma(e_{\sigma^{-1}}b)= (\sigma \cdot (e_{\sigma^{-1}}a))(\sigma \cdot (e_{\sigma^{-1}}b))\\
&\equal{\equref{4.1.0}}&
(\sigma\cdot e_{\sigma^{-1}})(\sigma\cdot a)(\sigma\cdot e_{\sigma^{-1}})(\sigma\cdot b)
=e_\sigma(\sigma\cdot a)e_\sigma(\sigma\cdot b)\\
&=& (\sigma\cdot 1_A)(\sigma\cdot a)(\sigma\cdot 1_A)(\sigma\cdot b)\equal{\equref{4.1.0}}
\sigma\cdot (ab)\equal{\equref{5.2.5}}\alpha_\sigma (e_{\sigma^{-1}}ab).
\end{eqnarray*}
We are left to prove that \equref{5.2.1} holds:
$$
\alpha_\sigma(e_{\sigma^{-1}}\alpha_\tau(e_{\tau^{-1}}a))
\equal{\equref{5.2.5}}\sigma\cdot(\tau\cdot a)\\
\equal{\equref{5.2.3}}e_\sigma((\sigma\tau)\cdot a)
\equal{\equref{5.2.5}}e_\sigma \alpha_{\sigma\tau}(e_{\tau^{-1}\sigma^{-1}}a).
$$
Conversely, assume that $G$ acts partially on $A$, and define an action of $kG$
on $A$ by extending
$\sigma\cdot a=\alpha_\sigma(e_{\sigma^{-1}}a)\in e_\sigma A$
linearly to $kG$. This defines a partial action of $kG$ on $A$, since
\begin{eqnarray*}
&&\hspace*{-2cm}
\sigma\cdot (ab)=\alpha_\sigma(e_{\sigma^{-1}}ab)
\equal{\equref{5.2.2}}
\alpha_\sigma(e_{\sigma^{-1}}a)\alpha_\sigma(e_{\sigma^{-1}}b)=(\sigma\cdot a)(\sigma\cdot b);\\
&&\hspace*{-2cm}
\sigma\cdot(\tau\cdot a)=\alpha_\sigma(e_{\sigma^{-1}}\alpha_\tau(e_{\tau^{-1}}a))
\equal{\equref{5.2.1}} e_\sigma\alpha_{\sigma\tau}(e_{\tau^{-1}\sigma^{-1}}a)
=e_{\sigma}((\sigma\tau)\cdot a);\\
&&\hspace*{-2cm} 1\cdot a=\alpha_1(e_1a)=\alpha_1(a)=a.
\end{eqnarray*}
It is easy to check that condition \equref{5.2.2b} establishes the bijectivity of the correspondence.
\end{proof}

\subsubsection*{Partial $U(L)$-actions}
The following result was kindly communicated to us by the referee.

\begin{proposition}\prlabel{5.2'}
Let $L$ be a Lie algebra over $k$ and $A$ a $k$-algebra. Let $H=U(L)$ be the universal enveloping algebra of $L$. Then every partial $H$-action on $A$ is an $H$-action on $A$ in the usual sense. 
\end{proposition}

\begin{proof}
One direction is clear. Conversely, suppose that $H=U(L)$ acts partially on $A$, i.e. conditions (\ref{eq:4.1.0},\ref{eq:4.1.1},\ref{eq:4.1.2}) are satisfied. For $x\in L$, condition (\ref{eq:4.1.0}) gives
\begin{eqnarray*}
&&\hspace{-1cm}
x\cdot (ab)=(x_{(1)}\cdot a)(x_{(2)}\cdot b)=(x\cdot a)(1\cdot b)+(1\cdot a)(x\cdot b)\equal{\equref{4.1.2}} x\cdot a + x\cdot b,
\end{eqnarray*}
for all $a,b\in A$. In particular, $x\cdot 1_A=x\cdot 1_A+x\cdot 1_A$, and thus $x\cdot 1_A=0$. \equref{4.1.1} then implies that, for all $x,y\in L$,
\begin{eqnarray*}
&&\hspace{-1cm}
x\cdot (y\cdot a)=(x_{(1)}\cdot 1_A)((x_{(2)}y)\cdot a)= (x\cdot 1_A)(y\cdot a)+(1\cdot 1_A)((xy)\cdot a)\\
&=&  (1\cdot 1_A)((xy)\cdot a)\equal{\equref{4.1.2}} (xy)\cdot a,
\end{eqnarray*}
hence \equref{4.2.2} follows. Also \equref{4.3.1} follows easily: $x\cdot 1_A=0=\epsilon(x)1_A$. We conclude that $A$ is an $H$-module algebra.
\end{proof}

\subsubsection*{A Frobenius property}
Let $i: R\to S$ be a ring homomorphism. Recall that $i$ is called Frobenius (or we say
that $S/R$ is Frobenius) if there exists
a Frobenius system $(\nu,e)$. This consists of an $R$-bimodule map $\nu: S\to R$
and an element $e=\sum e^1\ot_Re^2\in S\ot_R S$ such that
$se=es$, for all $s\in S$, and $\sum \nu(e^1)e^2=\sum e^1\nu(e^2)=1$.\\
A Hopf algebra $H$ over a commutative ring $k$ is Frobenius if and only if it is
finitely generated projective, and the space of integrals is free of rank one. If $H$
is Frobenius, then there exists a left integral $t\in H$ and a left integral $\varphi\in H^*$
such that
$
\lan \varphi, t\ran=1$.
The Frobenius system is $(\varphi, t_{(2)}\ot \ol{S}(t_{(1)}))$. In particular, we have
\begin{equation}\eqlabel{5.3.1a}
\lan\varphi, t_{(2)}\ran\ol{S}(t_{(1)})=t_{(2)}\lan\varphi,\ol{S}(t_{(1)})\ran=1_H.
\end{equation}
For a detailed discussion, we refer to the literature, see for example \cite[Sec. 3.2]{CMZ}.
If $t\in H$ is a left integral, then it is easy to prove that
\begin{equation}\eqlabel{5.3.2}
t_{(2)}\ot \ol{S}(t_{(1)})h=ht_{(2)}\ot \ol{S}(t_{(1)}),
\end{equation}
for all $h\in H$ (see \cite[Prop. 58]{CMZ} for a similar statement).\\
Assume that $A$ is a left $H$-module algebra, and that $H$ is Frobenius. Then
the ring homomorphism $A\to A\#H$ is also Frobenius (see \cite[Prop. 5.1]{CMZ}).
Similar properties hold for a module algebra over a weak Hopf algebra
and for an algebra with a partial group action (see \cite{CaenepeelDG05,CDG2}).
Our aim is now to prove such a statement for a partial module algebra over a
Frobenius Hopf algebra $H$. Assume that we have an action of $H$ on an algebra
$A$ satisfying (\ref{eq:4.1.0},\ref{eq:4.1.1},\ref{eq:4.1.2}). The smash product $A\# H$
has multiplication rule
$
(a\# h)(b\# g)=a(h_{(1)}\cdot b)\# h_{(2)}g$,
and $\ul{A\# H}$ is the subalgebra generated by the elements of the form
$(a\# h)1_A=a(h_{(1)}\cdot 1_A)\# h_{(2)}$.

\begin{proposition}\prlabel{5.3}
Let $H$ be a Frobenius Hopf algebra, let $t$ and $\varphi$ be as above, and take a left partial $H$-module algebra $A$. Suppose that $h\cdot 1_A$ is central in $A$, for every $h\in H$,
and that $t$ satisfies the following cocommutativity property:
\begin{equation}\eqlabel{5.3.4}
t_{(1)}\ot t_{(2)}\ot t_{(3)}\ot t_{(4)}=t_{(1)}\ot t_{(3)}\ot t_{(2)}\ot t_{(4)}.
\end{equation}
Then $\ul{A\# H}/ A$ is Frobenius, with Frobenius system
$(\ul{\nu}=(A\#\varphi)\circ \iota,
\ul{e}=(1_A\#t_{(2)})1_A\ot_A (1_A\# \ol{S}(t_{(1)}))1_A)$,
where $\iota: \ul{A\# H}\to A\# H$ is the inclusion map.
\end{proposition}

\begin{proof}
Applying $\Delta$ to the first tensor factor of \equref{5.3.4}, we see that
\begin{equation}\eqlabel{5.3.5}
t_{(1)}\ot t_{(2)}\ot t_{(3)}\ot t_{(4)}\ot t_{(5)}=t_{(1)}\ot t_{(2)}\ot t_{(4)}\ot t_{(3)} \ot t_{(5)}.
\end{equation}
For all $a\in A$ and $h\in H$, we have
\begin{eqnarray*}
&&\hspace{-1cm}(1_A\#t_{(2)})1_A \ot_A (1_A\# \ol{S}(t_{(1)}))1_A (a\#h)1_A\\
&=& (1_A\#t_{(2)}) \ot_A (1_A\# \ol{S}(t_{(1)}))(a\#h)1_A   \\
&=& (1_A\#t_{(2)}) \ot_A ((1_A\# \ol{S}(t_{(1)}))(a(h_{(1)}\cdot 1_A)\#h_{(2)}))1_A\\
&=& (1_A\#t_{(3)}) \ot_A (\ol{S}(t_{(2)})\cdot (a(h_{(1)}\cdot 1_A))\#\ol{S}(t_{(1)})h_{(2)})1_A   \\
&=& \bigl(t_{(3)}\cdot (\ol{S}(t_{(2)})\cdot (a(h_{(1)}\cdot 1_A))) \#t_{(4)}\bigr) \ot_A (1_A\# \ol{S}(t_{(1)})h_{(2)})1_A\\
&\equal{\equref{4.1.1}}& \bigl((t_{(3)}\cdot 1_A)((t_{(4)}\ol{S}(t_{(2)}))\cdot (a(h_{(1)}\cdot 1_A)))\#t_{(5)}\bigr)\ot_A (1_A\# \ol{S}(t_{(1)})h_{(2)})1_A \\
&\equal{\equref{5.3.5}}&\bigl( (t_{(4)}\cdot 1_A)((t_{(3)}\ol{S}(t_{(2)}))\cdot (a(h_{(1)}\cdot 1_A)))\#t_{(5)}\bigr)\ot_A (1_A\# \ol{S}(t_{(1)})h_{(2)})1_A\\
&=& \bigl((t_{(2)}\cdot 1_A)a(h_{(1)}\cdot 1_A) \#t_{(3)}\bigr) \ot_A (1_A\# \ol{S}(t_{(1)})h_{(2)})1_A\\
&=& \bigl(a(h_{(1)}\cdot 1_A) (t_{(2)}\cdot 1_A)\#t_{(3)}\bigr) \ot_A (1_A\# \ol{S}(t_{(1)})h_{(2)})1_A \\
&\equal{\equref{5.3.2}}& \bigl(a(h_{(1)}\cdot 1_A)((h_{(2)}t_{(2)})\cdot 1_A)\#h_{(3)}t_{(3)}\bigr) \ot_A (1_A\# \ol{S}(t_{(1)}))1_A \\
&=& (a(h_{(1)}\cdot 1_A)\#h_{(2)}t_{(2)})1_A \ot_A (1_A\# \ol{S}(t_{(1)}))1_A \\
&=& (a\#h)(1_A\#t_{(2)}) \ot_A (1_A\# \ol{S}(t_{(1)}))1_A \\
&=& (a\#h)1_A (1_A\#t_{(2)})1_A \ot_A (1_A\# \ol{S}(t_{(1)}))1_A.
\end{eqnarray*}
Using the fact that $\varphi$ is a left integral, we easily find that
$$\ul{\nu}((a\# h)1_A)=\lan \varphi, h_{(2)}\ran a (h_{(1)}\cdot 1_A)=
\lan \varphi, h\ran a(1_H\cdot 1_A)\equal{\equref{4.1.2}}\lan \varphi, h\ran a.$$
The left $A$-linearity of $\ul{\nu}$ is obvious, and the right $A$-linearity can be established
as follows:
\begin{eqnarray*}
&&\hspace*{-15mm}
 \ul{\nu}((a\# h)1_Ab)=\ul{\nu}((a\# h)b1_A)=\ul{\nu}((a(h_{(1)}\cdot b)\# h_{(2)})1_A)\\
 &=& \lan \varphi,h_{(2)}\ran a (h_{(1)}\cdot b)=\lan\varphi,h\ran a (1_H\cdot b)
 \equal{\equref{4.1.2}}\lan\varphi,h\ran ab=  \ul{\nu}((a\# h)1_A)b.
 \end{eqnarray*}
Finally,
\begin{eqnarray*}
&&\hspace*{-2cm} \ul{\nu} (1_A\#t_{(2)})1_A) ((1_A\# \ol{S}(t_{(1)}))1_A) = 
(\lan \varphi,t_{(2)}\ran 1_A\# \ol{S}(t_{(1)}))1_A\\
&\equal{\equref{5.3.1a}}& (1_A\# 1_H)1_A= 1_A\# 1_H;
\end{eqnarray*}
\begin{eqnarray*}
&&\hspace*{-2cm} (1_A\#t_{(2)}) \ul{\nu} ((1_A\# \ol{S}(t_{(1)}))1_A) 
=(1_A\# t_{(2)})\lan\varphi,\ol{S}(t_{(1)})\ran 1_A\\
&\equal{\equref{5.3.1a}}& (1_A\# 1_H)1_A= 1_A\# 1_H.
\end{eqnarray*}
\end{proof}

\begin{remark}
It follows from \equref{5.3.4} that $t$ is cocommutative. Obviously  \equref{5.3.4} is 
satisfied if $H$ is cocommutative.
\end{remark}

\section{Partial Hopf-Galois theory}\selabel{7}
Let $A$ be a right partial $H$-comodule algebra, and consider the corresponding $A$-coring
$\ul{\Cc}=\ul{A\ot H}$. The element $x=\rho(1_A)\in\ul{\Cc}$ is grouplike since
$\varepsilon(x)=\epsilon(1_{[1]})1_{[0]}=1_A$ (using \equref{2.6.1}) and
\begin{eqnarray*}
&&\hspace*{-15mm}
\Delta(x)=
(1_{[0]}\ot 1_{[1](1)})\ot_A (1_A\ot 1_{[1](2)})\\
&=&(1_{[0]}1_{[0']}\ot 1_{[1](1)}1_{[1']})\ot_A (1_A\ot 1_{[1](2)})\\
&\equal{\equref{2.2.1}}&
(1_{[0]}\ot 1_{[1]})\ot_A (1_A\ot 1_{[2]})
\equal{\equref{2.1.1}}
(1_{[0']}1_{[0]}\ot 1_{[1']}1_{[1]})\ot_A (1_A\ot 1_{[2]})
\\
&=&(1_{[0']}\ot 1_{[1']})\ot_A (1_{[0]}\ot 1_{[1]}) =x\ot_Ax.
\end{eqnarray*}
Then $A\in \Mm_A^H$.
Let $T=A^{{\rm co}H}=\{b\in A~|~\rho(b)=b\rho(1_A)\}.$
We have a morphism of corings
$\can: A\ot_TA\to \ul{A\ot H},~~\can(a\ot b)=ab_{[0]}\ot b_{[1]}$.
From \cite[Sec. 1 and Prop. 3.8]{Caenepeel03}
and \cite[Sec. 3]{W3}
we obtain immediately the following.

\begin{proposition}\colabel{7.2}
Let $A$ be a right partial $H$-comodule algebra. We have
a pair of adjoint functors $(F,G)$ between the categories $\Mm_T$ and
$\Mm_A^H$. $G=(-)^{{\rm co}H}$, and $F(N)=N\ot_TA$ with coaction
$\rho(n\ot_T a)=(n\ot_T a_{[0]})\ot a_{[1]}$. The following conditions are equivalent:\\
1) $\can$ is an isomorphism and $A$ is faithfully flat as a left $T$-module;\\
2) $(F,G)$ is an equivalence of categories and 
$A$ is flat as a left $T$-module;\\
3) $\ul{A\ot H}$ is flat as a left $A$-module, and $A$ is a projective generator of
$\Mm_A^H$.\\
In this situation, we will say that $A$ is a faithfully flat partial $H$-Galois extension of $T$.
\end{proposition}

From now on, we assume that $H$ is finitely generated and projective as a $k$-module,
with finite dual basis $\{(h_i,h_i^*)~|~i=1,\ldots,n\}$. Then $H$ is a $k$-progenerator: $H$ is a generator, because $\epsilon(1_H)=1$
(see for example \cite[I.1]{DI} for a discussion of (pro)generator modules). Suppose in addition that $A$ is finitely generated and projective as a left $T$-module.
Then $\can$ is an isomorphism if
and only if its left dual
${}^*\can: {}^*(\ul{A\ot H})\cong \ul{\#}(H,A)\to {}^*(A\ot_T A)\cong {}_T\End(A)^{\rm op}$
is an isomorphism. Viewed as a map $\ul{\#}(H,A)\to {}_T\End(A)^{\rm op}$, ${}^*\can$ is
given by the formula
${}^*\can(f)(a)=a_{[0]} f(a_{[1]})$. Composing ${}^*\can$ with the isomorphism
${}^*(\ul{A\ot H})^{\rm op}\cong \ul{A^{\rm op}\# H^{*{\rm cop}}}$ (see \prref{4.10}), we obtain
an $A^{\rm op}$-ring isomorphism $\theta: \ul{A^{\rm op}\# H^{*{\rm cop}}}\to {}_T\End(A)$.
We compute the map $\theta$ explicitly:
\begin{eqnarray*}
&&\hspace*{-15mm}
\theta((a\# h^*)\cdot 1_A)(b)=\theta(a\cdot (h^*_{(2)}\leftact 1_A)\# h^*_{(1)})(b)\\
&=& ({}^*{\rm can}\circ \alpha)((h^*_{(2)}\leftact 1_A)a\# h^*_{(1)})(b)=b_{[0]}\alpha((h^*_{(2)}\leftact 1_A)a\# h^*_{(1)})(b_{[1]})\\
&=& b_{[0]}(h^*_{(2)}\leftact 1_A)ah^*_{(1)}(b_{[1]})= b_{[0]}h^*_{(2)}(1_{[1]})1_{[0]}ah^*_{(1)}(b_{[1]})\\
&=& h^*(b_{[1]}1_{[1]})b_{[0]}1_{[0]}a \equal{\equref{2.1.1}}h^*(b_{[1]})b_{[0]}a.
\end{eqnarray*}
Recall (see \cite[Sec. 3]{CVW}) that we can associate a Morita context
$(T,{}^*\ul{\Cc},A,Q,\tau,\mu)$ to an $A$-coring $\ul{\Cc}$ with a fixed grouplike element $x$. We
will now compute this Morita context for $\ul{\Cc}=\ul{A\ot H}$ and 
$x=\rho(1_A)$, in the case where $A$ is a right partial $H$-comodule algebra.
First recall that $Q=\{q\in {}^*\ul{\Cc}~|~c_{(1)}q(c_{(2)})=q(c)x, \textrm{ for all $c\in \ul{\Cc}$}\}$. We first compute
$Q$ as a submodule of $\ul{\#}(H,A)$. $q\in \ul{\#}(H,A)$ satisfies the equation
\begin{equation}\eqlabel{7.3.1}
q(h)=1_{[0]} q(h1_{[1]}),
\end{equation}
for all $h\in H$. Let $\varphi$ be the map in ${}^*(\ul{A\ot H})$ corresponding to $q$. For $\gamma=
a1_{[0]}\ot h1_{[1]}\in \ul{A\ot H}$, we have that
$\Delta(\gamma)=(a1_{[0]}\ot h_{(1)}1_{[1]})\ot_A (1_A\ot h_{(2)}1_{[2]})$,
hence $q\in Q$ if and only if
$\varphi(\gamma)x=a1_{[0]} q(h1_{[1]})(1_{[0]}\ot 1_{[1]})=aq(h)1_{[0]}\ot 1_{[1]}$
equals
\begin{eqnarray*}
&&\hspace*{-2cm}
\gamma_{(1)}\varphi(\gamma_{(2)})
=(a1_{[0]}\ot h_{(1)}1_{[1]})q(h_{(2)}1_{[2]})\\
&=& a1_{[0]}q(h_{(2)}1_{[2]})_{[0]} \ot h_{(1)}1_{[1]}q(h_{(2)}1_{[2]})_{[1]} \\
&\equal{\equref{2.1.1}}& a(1_{[0]} q(h_{(2)}1_{[1]}))_{[0]} \ot h_{(1)}(1_{[0]}q(h_{(2)}1_{[1]}))_{[1]}\\
&\equal{\equref{7.3.1}}&aq(h_{(2)})_{[0]} \ot h_{(1)}q(h_{(2)})_{[1]},
\end{eqnarray*}
for all $a\in A$ and $h\in H$. We conclude that $Q$ is the submodule of $\#(H,A)$
consisting of the maps $q$ that satisfy \equref{7.3.1} and, for all $h\in H$,
\begin{equation}\eqlabel{7.3.2}
q(h_{(2)})_{[0]} \ot h_{(1)}q(h_{(2)})_{[1]}=q(h)1_{[0]}\ot 1_{[1]},
\end{equation}

Now we want to describe $Q$ as a submodule of $\ul{A^{\rm op}\#H^{*{\rm cop}}}\cong
\ul{\#}(H,A)^{\rm op}$. Take $\zeta=\sum_j a_j\# g_j^*\in \ul{A^{\rm op}\#H^{*{\rm cop}}}$
corresponding to $q\in \ul{\#}(H,A)^{\rm op}$. Then
\begin{equation}\eqlabel{7.3.3}
\sum_j a_j\# g_j^*=\sum_j (g_{j(2)}^*\leftact 1_A) a_j\# g_{j(1)}^*
\end{equation}
where $\leftact$ is the map from the left partial $H^{*{\rm cop}}$-action on $A^{\rm op}$, corresponding to
$\rho$, cf. \thref{4.8}. Then $\zeta\in Q$ if and only if
$\sum_j a_{j[0]} g_j^*(h_{(2)}) \ot h_{(1)}a_{j[1]}=
\sum_j a_jg_j^*(h) 1_{[0]} \ot 1_{[1]}$,
for all $h\in H$. This is equivalent to stating that
\begin{eqnarray*}
&&\hspace*{-2cm}
\sum_j a_{j[0]} g_j^*(h_{(2)}) h^*(h_{(1)}a_{j[1]})=
\sum_j a_{j[0]} g_j^*(h_{(2)}) h^*_{(2)}(h_{(1)})h^*_{(1)}(a_{j[1]})\\
&\equal{\equref{4.6.2}}& \sum_j h^*_{(1)}\leftact a_j g_j^*(h_{(2)}) h^*_{(2)}(h_{(1)})=
\sum_j h^*_{(1)}\leftact a_j (h^*_{(2)}*g_j^*)(h)
\end{eqnarray*}
is equal to
\begin{eqnarray*}
&&\hspace*{-2cm}
\sum_j a_jg_j^*(h) 1_{[0]} h^*(1_{[1]})
= \sum_{j} a_j1_{[0]}g_j^*(h) h^*_{(2)}(1_{[1]})h^*_{(1)}(1_H)\\
&\equal{\equref{4.6.2}}&
\sum_j a_j (h^*_{(2)}\leftact 1_A) g^*_j(h) h^*_{(1)}(1_H),
\end{eqnarray*}
for all $h\in H$ and $h^*\in H^*$. We conclude that $Q$ consists of the elements
$\zeta=\sum_j a_j\# g_j^*\in {A^{\rm op}\#H^{*{\rm cop}}}$ satisfying \equref{7.3.3} and
\begin{equation}\eqlabel{7.3.4}
\sum_j h^*_{(1)}\leftact a_j\# (h^*_{(2)}*g_j^*)=\sum_j a_j (h^*_{(2)}\leftact 1_A)\# h^*_{(1)}(1_H)g^*_j,
\end{equation}
for all $h^*\in H^{*{\rm cop}}$.
It follows from \cite[Lemma 3.1]{CVW} that $Q$ is a $(\ul{\#}(H,A),T)$-bimodule.
The left action is given by the multiplication in $\ul{\#}(H,A)$. We also know
(cf. \cite[Prop. 2.2]{CVW}) that $A$ is a $(T,\ul{\#}(H,A))$-bimodule, with right $\ul{\#}(H,A)$-action
given by the formula $a\cdot f=a_{[0]} f(a_{[1]})$. We have well-defined maps
$$\tau: A\ot_{\ul{\#}(H,A)} Q\to T,~~\tau(a\ot q)=a\cdot q=a_{[0]} q(a_{[1]});$$
$$\mu: Q\ot_TA\to \ul{\#}(H,A),~~\mu(q\ot a)(h)=q(h)a.$$
$(T,\ul{\#}(H,A),A,Q,\tau,\mu)$ is a Morita context. The map $\tau$ is surjective if and only
if there exists $q\in Q$ such that $q(1_H)=1_A$. 

\begin{theorem}\thlabel{7.3}
Let $A$ be a right partial $H$-comodule algebra, and assume that $H$ is finitely generated
and projective
as a $k$-module. Then the following assertions are equivalent.
\begin{enumerate}
\item $A$ is a faithfully flat partial $H$-Galois extension of $T$;
\item $\theta$ is an isomorphism and $A$ is a left $T$-progenerator;
\item the Morita context $(T,\ul{\#}(H,A),A,Q,\tau,\mu)$ is strict;
\item $(F,G)$ is an equivalence of categories.
\end{enumerate}
\end{theorem}

\begin{proof}
Since $H$ is finitely generated projective as a $k$-module, $A\ot H$ is finitely
generated and projective as a left $A$-module. Being a direct factor of
$A\ot H$, $\ul{A\ot H}$  is also finitely
generated and projective as a left $A$-module. $\ul{A\ot H}$ is a left $A$-generator
since
$\varepsilon(1_{[0]}\ot 1_{[1]})=1_{[0]}\epsilon(1_{[1]})
= 1_A$ (use \equref{2.2.3}).
It follows that $\ul{A\ot H}$ is a left $A$-progenerator, and the result then follows
immediately from \cite[Theorem 4.7]{Caenepeel03}.
\end{proof}

Assume that $H$ is Frobenius. Then $H^{*{\rm cop}}$ is also Frobenius. Assume,
moreover, that
$a1_{[0]}\ot 1_{[1]}=1_{[0]}a\ot 1_{[1]}$ and
$\lan \varphi,hgh'g'\ran= \lan \varphi,hh'gg'\ran$,
for all $a\in A$ and $h,g,h',g'\in H$. Then it follows from \prref{5.3} that
$\ul{A^{\rm op}\#H^{*{\rm cop}}}/A^{\rm op}$ is Frobenius. The Morita context 
$(T,\ul{\#}(H,A)\cong \ul{A^{\rm op}\#H^{*{\rm cop}}},A,Q,\tau,\mu)$ is the Morita
context associated to the $A^{\rm op}$-ring\\
 $\ul{A^{\rm op}\#H^{*{\rm cop}}}$,
see \cite[Theorem 3.5]{CVW}. It follows from \cite[Theorem 2.7]{CVW}
that $Q\cong A$ as $k$-modules. So we conclude that the Morita context is of the
form $(T, \ul{A^{\rm op}\#H^{*{\rm cop}}},A,A,\tau,\mu)$.

\end{document}